\newif\ifAcm

\Acmfalse 	

\ifAcm
\documentclass[acmsmall]{acmart}
\usepackage{hyperref}
\else
\documentclass{article}

\usepackage[utf8]{inputenc}
\usepackage[multiple]{footmisc}
\usepackage[pdftex,colorlinks,urlcolor=blue]{hyperref}
\fi

\usepackage[linesnumbered,ruled,vlined]{algorithm2e}
\usepackage[normalem]{ulem}
\usepackage[leqno]{amsmath}
\usepackage{amsthm}
\usepackage{amssymb}
\usepackage{todonotes}
\usepackage{mathtools}
\usepackage{amsmath}
\usepackage{enumerate}
\usepackage{pdflscape}
\usepackage{algorithmic}
\usepackage[font=footnotesize,labelfont=bf]{caption}
\usepackage[font=scriptsize,labelfont=bf]{subcaption}

\SetKwRepeat{Do}{do}{while}

\def\e{\mathbf{e}}
\def\b{\mathbf{b}}
\def\c{\mathbf{c}}

\def\x{\mathbf{x}}
\def\z{\mathbf{z}}

\def\vecgamma{\mathbf{\gamma}}
\def\st{\mathrm{s.t.}}
\def\Sy{\mathcal{S}}
\def\prob{\bf}
\def\conv{\mathrm{conv}}

\newcommand{\expedis}{{EXPEDIS}}

\newcommand{\psd}{\succcurlyeq 0}

\newcommand{\unaryminus}{\scalebox{0.7}[1.0]{\( - \)}}
\newcommand{\hs}{\{\unaryminus 1,1\}^{n}}

\DeclarePairedDelimiter{\card}{\lvert}{\rvert}

\newcommand\R{{\mathbb R}}

\newcommand\prods[2]{\langle #1, #2 \rangle}

\DeclareMathOperator{\diag}{diag}
\DeclareMathOperator{\Diag}{Diag}
\DeclareMathOperator{\OPT}{\mathrm{OPT}}
\DeclareMathOperator{\LB}{\mathrm{LB}}
\DeclareMathOperator{\rank}{rk}

\DeclareMathOperator{\SU}{{\tt \mbox{speed-up} }_{{\tt inst}}}
\DeclareMathOperator{\diff}{\mathrm{diff}}

\newtheorem{theorem}{Theorem}

\newcommand{\rew}{\textcolor{black}} 
\newcommand{\reww}{\color{black}}

\begin{document}
	
	\ifAcm
\else
	
	\title{BiqBin: a parallel branch-and-bound solver for binary quadratic problems with linear constraints}
	
	\author{Nicol{\`o} Gusmeroli\thanks{Alpen-Adria-Universit\"{a}t Klagenfurt, Austria. \newline Emails: \{nicolo.gusmeroli,melanie.siebenhofer,angelika.wiegele\}@aau.at}
		\and{Timotej Hrga}\thanks{University of Ljubljana, Faculty of Mechanical Engineering, Ljubljana, Slovenia. \newline Emails: \{timotej.hrga,janez.povh\}@lecad.fs.uni-lj.si}
		\and{Borut Lu\v{z}ar}\thanks{Faculty of Information Studies, Novo mesto, Slovenia.	Email: borut.luzar@gmail.com}
		\and{Janez Povh}\footnotemark[2]
		\and{Melanie Siebenhofer}\footnotemark[1]
		\and{Angelika Wiegele}\footnotemark[1]}
	\date{\today}
	
	\maketitle
	
	\fi
	
	\begin{abstract}
		We present BiqBin, an exact solver for linearly constrained binary quadratic problems. Our approach is based on an exact penalty method to first efficiently transform the original problem into an instance of Max-Cut, and then to solve the Max-Cut problem by a branch-and-bound algorithm. All the main ingredients are carefully developed using new semidefinite programming relaxations obtained by strengthening the existing relaxations  with a set of hypermetric inequalities, applying the bundle method as the bounding routine and using new strategies for exploring the branch-and-bound tree.
		
		Furthermore, an efficient C implementation of a sequential and a parallel branch-and-bound algorithm is presented.
		The latter is based on a load coordinator-worker scheme using MPI for multi-node parallelization and is evaluated on a high-performance computer.
		
		The new solver is benchmarked against BiqCrunch, GUROBI, and SCIP on four families of (linearly constrained) binary quadratic problems.
		Numerical results demonstrate that BiqBin is a highly competitive solver. The serial version outperforms the other three solvers on the majority of the benchmark instances.
		We also evaluate the parallel solver and show that it has good scaling properties. The general audience can use it as an on-line service available at \url{http://www.biqbin.eu}.
		
	\end{abstract}

\section{Introduction}
	
\subsection{Motivation}\label{sec:mot}
With the advent of data driven economy,
a close cooperation between data science and mathematical optimization became a crucial driver
for developing (i) new data science methods that are capable to reveal new knowledge hidden in the data
and (ii) new algorithms for implementing these data science methods.

Studying how to group data instances according to their inner similarity,
i.e., data clustering analysis~\cite[Ch.10.3]{GaWiHa:13},
is a traditional and one of the most studied problems in data science.
When we know that the data should be decomposed into a fixed number of groups, say $K$,
and we want to find these $K$ groups, we have the $K$-clustering problem,
which is NP-hard~\cite{GaJo:79}.

Another interesting problem from data science, similar to the $K$-clustering problem, is the problem where
we want to find a set of $k$ vertices in a simple undirected weighted graph $G$ such that the sum of the weights on the edges connecting these $k$ vertices is maximum.
	This problem is called the densest $k$-subgraph problem \cite{billionnet2005different}
and can be formulated as an optimization problem with a quadratic objective function, one linear constraint,
and binary ($0/1$) decision variables:
\begin{align}
\max~~ & \frac{1}{2}\x^\top  W \x\nonumber\\
\st~~ & \e^\top \x = k\tag{\prob{D$k$S}} \label{eqn:DkS}\\
& \x \in \{0,1\}^{n}.\nonumber
\end{align}
Here, $W$ is the weighted adjacency matrix of the underlying graph and $\e$ is the all one vector.

The densest $k$-subgraph problem can be seen as a generalization of the Max-Clique problem~\cite{wu2015review}
and also as a special case of the quadratic knapsack problem~\cite{pisinger2007quadratic}.
Even though it is easy to understand and it appears to be simple at first sight, it is not.
In fact, it is one of the NP-hard problems for which even the approximability is not well understood.
There is a huge gap between the best approximation algorithm and known inapproximability results~\cite{bhaskara2012polynomial}.

The problems mentioned above are all special instances of binary quadratic problems with linear constraints, which are formally defined in Section~\ref{sec:mc}.
In practice, we typically solve such problems only approximately using heuristic algorithms.
For example, in~\cite{das2009metaheuristic,wiwie2015comparing,XuWu:05} one can find a vast amount of such heuristics for the case of data clustering.
However, to evaluate the performance of the heuristic algorithms, we still need the ground truth,
i.e., the optimum solutions of the original problems.
Therefore, solving (to optimality) constrained binary quadratic problems, like~\eqref{eqn:DkS}, on as large as possible problem instances is highly needed.

Many other optimization problems with clear real-life applications can also be represented in a similar way as a
non-convex optimization problem in binary variables subject to linear constraints,
e.g., the quadratic assignment problem, the stable set problem, and the graph coloring problem (see, e.g.,~\cite{schrijver2003combinatorial} for definitions).
Again, appropriate (meta)heuristic algorithms are used to approximately solve instances of big size,
while optimum solutions for test instances of small or medium size are still needed to evaluate these heuristics.

\subsection{Notation}

We will use the following notation:
$\Sy_n$ denotes the space of symmetric $n\times n$ matrices, $\R^n$ is
the $n$-dimensional space of all $n$-tuples of real numbers and
$\R^{m\times n}$ is the space of all $m\times n$ real matrices.
By $\e_n$ we denote the vector of length $n$ with all entries equal to one.
Usually, its dimension is clear from the context, so we write only $\e$.
Similarly, $\mathbf{0}$ denotes the zero vector or the zero matrix.
Given $\x \in \R^n$, $\Diag(\x)$ is the $n \times n$ diagonal matrix
with $\x$ on its diagonal;
and $\diag(X)$ is the vector with the diagonal elements of matrix $X$.
By $\rank(X)$ we denote the rank of matrix $X$.

For an optimization problem $\prob{P}$, we refer to its optimum value by  $\OPT_{\prob{P}}$.
For the SDP relaxations of  \eqref{eqn:MC}
 we use slightly different notations for their optimum values.
For example, the optimum value of the SDP
relaxation \eqref{basic_SDP} is denoted by $\OPT_{\prob{SDP}}$.

\subsection{Linearly constrained binary quadratic problem and the Max-Cut problem}	
\label{sec:mc}

The central problem that we consider in the paper is the linearly constrained binary quadratic problem (BQP), which can be formulated as follows:
\begin{align}
	\label{eqn:BQP}  \tag{\prob{BQP}}
	\begin{split}
	\min~~  &\x^\top F \x+ \c^\top \x \\
	\st~~   &A\x = \b,~~ \x \in \{0,1\}^{n},
	\end{split}
\end{align}
with given data $F \in \Sy_n$, $c\in \R^{n}$, $A\in \R^{m\times n}$, and $\b\in \R^{m}$.
We assume to have only linear equality constraints, in case of inequalities we may add a slack variable
and decompose it into a weighted sum of boolean variables.
This problem encompasses also the densest $k$-subgraph problem and is thus an NP-hard problem.

The mathematical optimization community is interested in BQP as a problem \emph{per se}.
In this context, the main challenge is developing new algorithms for it by
(i) exploring and exploiting new properties of the problem,
(ii) using new results from other (non-optimization) areas (such as algebraic geometry),
(iii) making new combinations of existing algorithms with best practical or theoretical performance,
and (iv) exploiting best available high-performance hardware and software.

Lasserre~\cite{lasserre2016max} has proved that any instance of \eqref{eqn:BQP}
can be transformed into an instance of the Max-Cut problem,
which is an NP-hard optimization problem~\cite{GaJo:79,Ka:72} on graphs.
This problem is among the most studied combinatorial optimization problems, it has connections to various fields of discrete mathematics and models a wide range of applications.
The transformation is based on an exact penalty approach, that was further explored and advanced in a recent paper of two co-authors of this paper~\cite{GuWi:19}.

The Max-Cut problem can be defined as follows.
Suppose a weighted undirected graph $G=(V,E)$ is given,
where $V$ is the set of vertices, $E$ is the set of edges, and each edge $e\in E$ has the weight $w_e\in \R$.
The Max-Cut problem asks to find a partition of $V$ into two parts $(S, V\backslash S)$ such that the sum of the weights
of the edges having one endpoint in $S$ and the other one in $V\backslash S$ is maximized.

If $W=(w_{ij})$ is the weighted adjacency matrix of $G$, i.e., $w_{ij} = w_{ji} = w_e$ for $e=\{i,j\} \in E$,
and the Laplacian matrix of the graph associated with $W$ is
\begin{equation*}
	L = \Diag(W\e) - W,
\end{equation*}
then computing the maximum cut amounts to solving the following binary quadratic problem in variables from $\{\unaryminus 1,1\}$:
\begin{equation*}
	\label{eqn:MC}  \tag{\prob{Max-Cut}}
	\OPT_{\prob{Max-Cut}}~=~\max \left\{ \frac{1}{4}\x^\top L\x \mid \x\in
	\{\unaryminus 1,1\}^{\card{V}} \right\}.
\end{equation*}
However, it is straightforward to transform it into a $0/1$ problem by using the following simple linear substitution $\x = 2\z - \e$ and by observing that
$$
	\x^\top L\x = 4\z^\top L\z -4\z^\top L\e +\e^\top L\e = 4\z^\top L\z,
$$
since $L\e =0 $.
Hence, we obtain the equivalent formulation of the Max-Cut, which has a structure of~\eqref{eqn:BQP} (with no linear constraints):
\begin{equation*}
\label{eqn:MC2}
\OPT_{\prob{Max-Cut}}~=~\max \left\{ \z^\top L\z \mid \z \in \{0,1\}^{\card{V}} \right\}.
\end{equation*}

The reason why the community is interested in reformulating \eqref{eqn:BQP} into \eqref{eqn:MC} is related to the fact that for the latter problem,
there exists a wide range of approximate and exact methods and solvers and we want to employ them in solving the former problem.
We are particularly interested in the methods computing global optima of both problems.

\subsection{Our contribution}

The main contribution of this paper is the BiqBin solver for BQP,
which outperforms the existing solvers on some special cases of BQP.
The core of the BiqBin solver is the exact penalty algorithm \expedis~\cite{GuWi:19},
meaning that we first transform every instance of~\eqref{eqn:BQP} to a corresponding instance of~\eqref{eqn:MC} and then solve the resulting instance of~\eqref{eqn:MC}.
BiqBin is coded in C and non-trivially improves existing solvers by
\begin{itemize}
	\item{} introducing a strengthened bounding routine based on hypermetric inequalities;
	\item{} implementing a parallel branch-and-bound algorithm to solve \eqref{eqn:MC} instances using a Message Passing Interface library (MPI);
	\item{} providing a web-based BiqBin service\footnote{\url{http://www.biqbin.eu}}
	enabling researchers to submit their instances of~\eqref{eqn:BQP} to one of
	the Slovenian Tier-2 supercomputers to be solved.
\end{itemize}
Additionally, we demonstrate practical efficiency of BiqBin by providing an extensive benchmarking with
BiqCrunch~\cite{krislock2017biqcrunch}, GUROBI~\cite{gurobi} and SCIP~\cite{achterberg2009scip}
on the list of four special cases of BQP, including the Max-Cut problem, the unconstrained binary quadratic problem, the densest $k$-subgraph problem and
randomly generated binary quadratic problems with linear constraints.

We can observe that BiqBin is performing very well on the instances with a small number of linear constraints,
while on the instances with many linear constraints,
the solvers which directly exploit the structure of the constraints are very competitive.
Additionally, as we show on the benchmark instances, BiqBin is also highly scalable.

\section{Related work}
Starting in the previous century, tackling the Max-Cut problem computationally attracted attention of many researchers.
A number of ideas for solving the problem have been proposed in the literature;
fast heuristic algorithms on one side and exact algorithms on the other side have been developed over the last decades.
We mention here only a few concepts of such exact methods.
A first success on solving \eqref{eqn:MC} instances to optimality appeared in the eighties,
when linear programming based methods have been implemented~\cite{BaJuRe:89} and further developed;
in particular, in the context of solving problems arising in physics~\cite{LiJuReRi:04}.
These methods are in specifically successful in the cases when the underlying graphs are sparse.

Other methods use a preprocessing phase where they try to fix some variables (based on the gradient)~\cite{PaRo:90,PaRo:90b}
or to construct a convex problem having the same optimal solution. For such convex problems efficient solvers exist;
note, however, that most of these are commercial solvers, e.g., CPLEX~\cite{cplex}.

At the beginning of this century,
methods based on semidefinite programming and the branch-and-bound (B\&B) algorithm have been proposed~\cite{HeRe:98,KrMaRo:14,ReRiWi:10,krislock2017biqcrunch}.
One of them, the BiqMac solver~\cite{ReRiWi:10}, was the starting point for our research,
since it is one of the best performing solvers for \eqref{eqn:MC} instances
and can be employed to solve \eqref{eqn:BQP} instances
as follows from the Lassere's result~\cite{lasserre2016max}.
Another one, the BiqCrunch~\cite{krislock2017biqcrunch} solver, we use to benchmark our results.
			
By setting $F=0$ in~\eqref{eqn:BQP}, we obtain a linear optimization problem with binary variables.
Such problems have been investigated for many decades, leading to enormous progress in their practical solutions.
According to Hans Mittelmann's web page\footnote{\url{http://plato.asu.edu/ftp/milp.html}},
the state-of-the-art solvers (the web page contains results for the solvers {CBC, GLPK, LP\_SOLVE, MATLAB,  SAS-OR}, and {(F)SCIP})
can solve such problems to optimality within few hours if the number of binary variables is up to 50.000.
Moreover, they can also handle problems with several hundred thousand binary variables if sufficient structure is provided
(see also the list of ``easy'' cases on the MIPLIB benchmark page~\cite{miplib2017}).
	
In the case of a (non-convex) quadratic objective function,
the state-of-the-art optimization techniques perform weaker than those designed for linear problems.
There are many solvers that can solve \eqref{eqn:BQP} instances:
Mittelmann included in his decision tree for optimization software\footnote{\url{http://plato.asu.edu/sub/nlores.html\#QP-problem}} several solvers
which use various optimization techniques to compute optimum solutions,
such as the branch-and-cut algorithm, branch-and-bound, lift-and-project, convex reformulation combined with some first and second order methods, etc.
He performed an extensive benchmarking of several solvers ({BARON, (F)SCIP, ANTIGONE, MINOTAUR, OCTERACT, GUROBI})
on instances from QPLIB\footnote{\url{http://qplib.zib.de/index.html}}\footnote{\url{http://plato.asu.edu/ftp/qplib.html}}~\cite{furini2019qplib}.
Mittelmann showed that these solvers can solve within 1 hour between 5\% ({MINOTAUR}) and 63\% ({GUROBI}) of the 128 benchmarking instances.
These benchmarking problems have small to medium size: they mostly contain up to a few hundred binary variables with the largest instance having 8904 binary variables.
Large instances typically share some important structural properties, which make them solvable at least for some solvers.

The solvers mentioned above are available under various software licences.
Most of them are freely available to all researchers for academic purposes upon registration and verification.
If one does not have a strong local machine or does not want to bother with local installations,
they can submit the problem instance to the web portals where some of the solvers are installed;
such services are available, e.g., for NEOS~\cite{dolan2001neos}, BiqMac\footnote{\url{http://biqmac.aau.at/}},
and BiqCrunch\footnote{\url{https://biqcrunch.lipn.univ-paris13.fr/BiqCrunch/solver/}}.
While NEOS partially runs on fast supercomputers,
BiqMac and BiqCrunch run on a small cluster and a strong single machine with multicore processors, respectively.
The latter two solvers utilize B\&B algorithms, but they do not perform any parallelization
\rew{(a first trial to parallelize BiqMac was done in \cite{Kalamar_et_al_2017}).}
They are available on-line (more precisely, the users can submit their instance of the problem on-line) but they are running on a hardware,
which cannot be compared to the state-of-the-art supercomputers.

\section{Semidefinite programming relaxations for the Max-Cut}
\label{sec:sdp}

In the subsequent sections, we make use of tools from semidefinite programming.
In order to make this paper self-contained, we recall here some definitions and algorithms.

As mentioned in Section~\ref{sec:mc}, the Max-Cut problem can be formulated as
\begin{equation*}
	\max \left\{ \frac{1}{4}\x^\top L\x \mid \x\in \{\unaryminus 1,1\}^{\card{V}} \right\}.
\end{equation*}
Observe that for any $\x \in \{\unaryminus 1,1\}^{\card{V}}$, the matrix $X = \x\x ^\top $ is positive semidefinite and its diagonal is equal to the vector of all ones.
Using this transformation and the property $\x^\top L\x = \prods{L}{ \x\x^\top}$, we can re-write~\eqref{eqn:MC} as
\begin{equation*}
\max \left\{\frac{1}{4}\prods{L}{X} \mid  \diag(X) = \e, \ X \psd, \ \rank(X) = 1 \right\}.
\end{equation*}

By dropping the rank-one constraint, we obtain the basic SDP relaxation
\begin{equation}
	\label{basic_SDP}
	\tag{$\mathbf{MC}_{\prob{SDP}}$}
	\OPT_{\prob{SDP}}~=~\max \left\{\frac{1}{4}\prods{L}{X} \mid \ X \psd, \ \diag(X) = \e \right\}.
\end{equation}
It is well-known that the bound $\OPT_{\prob{SDP}}$ is not strong enough to be successfully used within the B\&B framework
even for solving the Max-Cut problem to optimality on graphs with only 50~nodes.
We overcome this problem
by adding additional equality or inequality constraints, known as cutting planes,
to strengthen the bound and consequently decrease the size of the B\&B tree.

Similarly as in BiqMac and BiqCrunch, we use triangle inequalities.
Furthermore, we strengthen the bound by using higher order $k$-gonal inequalities, which belong to the family of hypermetric inequalities~\cite{DeLa:94}.
They can be introduced as follows.
Suppose that $\b$ is an integer vector for which $\e^\top \b$ is odd. This implies that
$
\vert \x^\top \b \vert \ge 1 \textrm{ for all } \x \in \{\unaryminus 1,1\}^n
$
and therefore
$
\prods{\b\b^\top}{\x\x^\top} \ge 1.
$
The \emph{hypermetric inequalities} are the following set of linear inequalities,
which can be applied to any symmetric matrix $X$ of order $n$,
and are valid for any $X$ from the convex hull of rank-one matrices $\x\x^\top$, for $\x\in \{\unaryminus 1,1\}^n$:
\begin{equation*}
\{\prods{\b\b^\top}{X} \ge 1\ \mid \ \e^\top \b \textrm{ odd}, \ \b \textrm{ integer}\}.
\end{equation*}

In this paper, we consider the subclasses of hypermetric inequalities generated by choosing $\b$ with $b_i \in \{\unaryminus 1,0,1\}$
and by fixing the number of non-zero entries in $\b$ to 3, 5, or 7.
In this cases, we obtain triangle, pentagonal, and heptagonal inequalities, respectively.
The latter two are also called $5$-clique and $7$-clique inequalities, respectively.
\rew{
More specifically, the triangle inequalities are defined as
\begin{align*}
-X_{ij}-X_{ik}-X_{jk} \le 1, \\
-X_{ij}+X_{ik}+X_{jk} \le 1, \\
X_{ij}-X_{ik}+X_{jk} \le 1, \\
X_{ij}+X_{ik}-X_{jk} \le 1,
\end{align*}
$\forall$ distinct $i$, $j$, $k$.
The first inequality is actually $b^\top Xb = X_{ii} + X_{jj} + X_{kk} + 2\left(X_{ij} + X_{ik}+ X_{jk}\right) \ge 1$,
for the case $\b_i = \b_j = \b_k = 1$, by using the implicit constraint $\diag(X) = \e$.
The other three inequalities  from above, for selected $i,j,k$, follow by considering all the other combinations of signs of entries of $\b_i,\b_j,\b_k$ (trivially, $\b$ and $-\b$ yield the same inequality).
}
By introducing appropriate linear operators on the vector space of symmetric matrices
and by using the fact that the matrices under consideration have diagonal entries equal to one,
we write such inequalities as $\mathcal{A}_3(X)\le \e$, $\mathcal{A}_5(X)\le \e$ and $\mathcal{A}_7(X)\le \e$, respectively.

Similarly, we denote by $\mathcal{A}_{\prob{HYP}}(X) \le \e$ a set containing triangle, pentagonal, and heptagonal inequalities
and call it (by a slight abuse of notation) hypermetric inequalities.
We infer the following strengthening of \eqref{basic_SDP}:
\begin{equation}\label{hyper_SDP}\small
\tag{$\mathbf{MC}_{\prob{HYP}}$}
	\OPT_{\prob{HYP}}~=~\max \left\{\prods{L}{X} \mid \ X \psd, \ \diag(X) = \e, \ \mathcal{A}_{\prob{HYP}}(X)\le \e \right\}.
\end{equation}

We now describe the routine for separating the proposed cutting planes.
There are $4{n \choose 3}$ triangle inequalities;
for problem sizes that we are interested in,
we can enumerate them and identify the most violated ones.

Due to a large number of higher $k$-gonal inequalities,
separation of pentagonal and heptagonal inequalities is done heuristically.
Let $\e = (1, 1, 1, 1, 1)^\top$ and define $H_1 = \e\e^\top$.
Suppose we are searching for the pentagonal inequality
(of the type where all nonzero entries of $\b$ are ones)
with large violation, i.e., for a given matrix $X$,
we are looking for a 5-permutation $p$ of $n$ vertices such that the value
$
\prods{H_1}{X(p,p)}
$
is minimal.
By $X(p,p)$ we denote the submatrix obtained by taking rows and columns contained in the permutation $p$.
These numbers represent the indices of nonzero entries of $\b$ determining the pentagonal inequality.
Let $H$ be an $n \times n$ matrix having $H_1$ as the leading principal submatrix of order 5 and all other elements are set to zero.
Then the problem can be reformulated as a quadratic assignment problem of the form
\begin{align*}
\min \quad& \prods{H}{PXP^\top }\\
\st \quad& P \in \Pi,
\end{align*}
where $\Pi$ is the set of all $n \times n$ permutation matrices.
This problem is approximately solved by using simulated annealing to obtain a pentagonal inequality with potentially large violation.
By replacing the matrix $H_1$ with rank-one matrices
$H_2 = \widehat{\e}\widehat{\e}^\top$ or $H_3 = \widetilde{\e}\widetilde{\e}^\top $,
where $\widehat{\e} = (\unaryminus 1, 1, 1, 1, 1)$ and
$\widetilde{\e} = (\unaryminus 1, \unaryminus 1, 1, 1, 1)$,
different types of pentagonal inequalities are found.
The same idea is applied for separating \rew{strongly} violated heptagonal inequalities.

To sum up, in BiqBin, we first iteratively identify a subset of the promising cutting planes.
Using a current approximate solution $X$,
we find the \rew{promising} hypermetric inequalities $\mathcal{A}'_{\prob{HYP}}(X)\le \e $ and solve \eqref{hyper_SDP} by using only these inequalities.
Specifically, during the separation routine, we add $10\cdot n$ triangle inequalities, $300$ pentagonal inequalities and $200$ heptagonal inequalities.

Note that for convenience, we sometimes use instead of the cut vector
$\x\in \{\unaryminus 1,1\}^n$ the vector
$\reww{\begin{pmatrix} 1 \\ \x\end{pmatrix} \in \{\unaryminus 1,1\}^{n+1}}$
to derive the SDP relaxation.
This has the advantage that the values of the vector $\x$ are given
in the \rew{first} row and column of the matrix
$$
\reww{\begin{pmatrix} 1 \\ \x \end{pmatrix}
\begin{pmatrix} 1\\ \x \end{pmatrix}^\top.}
$$

\section{From binary quadratic problems with linear constraints to the Max-Cut problem}
\label{sec:exact_penalty}

In this section, we recall how to transform binary quadratic problems with linear constraints
into the Max-Cut problem.

Exact penalty methods for solving constrained optimization problems construct a function,
for which the (unconstrained) minimizers are also optimal solutions of the constrained problem,
see~\cite{di1994exact} for an overview of classical results on this topic.
Gusmeroli and Wiegele~\cite{GuWi:19} introduced an exact penalty algorithm {\em over discrete sets} called \expedis. 
Their work follows and improves the idea of Lasserre, see~\cite{lasserre2016max},
and reformulates a linearly constrained binary quadratic problem as a Max-Cut instance.

The input of the {\expedis} algorithm is an instance of \eqref{eqn:BQP}.
We consider its version with the binary variables being from $\{\unaryminus1,1\}$.

In order to simplify notation, we define the sets of feasible and infeasible binary vectors as $\Delta$ and $\Delta^c$, respectively;
i.e.,
$$
	\Delta = \left\{\x \in \hs \mid A\x=\b \right\} \quad \textrm{and} \quad \Delta^c = \hs \setminus \Delta.
$$

Given a sufficiently large penalty parameter, denoted $\sigma$, we add a quadratic term to the objective function and obtain
\[
	h(\x) = \x^\top F\x + \c^\top \x + \sigma \|A\x-\b\|^{2}.
\]
By defining the matrix
\[
	Q = \begin{bmatrix}
	\sigma b^\top b & \left( c - 2\sigma A^\top b \right)^\top/2 \\ \left(c - 2\sigma A^\top b\right)/2 & F + \sigma A^\top A
	\end{bmatrix}
\]
the function $h(\x)$ can alternatively be written as $\bar{\x}^\top Q \bar{\x}$, in which case
we consider the following unconstrained binary optimization problem
\begin{align*}
	\begin{split}
	h^* = \min~~ & \bar{\x}^\top Q \bar{\x}  \\
	\st~~ & \bar{\x} \in \left\{ \unaryminus 1,1 \right\}^{n+1} \\
	& \reww{\bar{\x}_0 = 1},
	\end{split}
\end{align*}
which is a Max-Cut problem on a graph with $n+1$ vertices\rew{, i.e.,
	$\lvert V \rvert = \{0,1, \dots, n\}$};
see~\cite{GuWi:19} for more details.

\begin{theorem}\cite[Theorem~2]{GuWi:19} 
	\label{thm:exp}
        Consider an instance of \eqref{eqn:BQP} with optimal value $f^*$. Let $\rho$ be a threshold parameter and let $\sigma$ be a penalty parameter such that
	\begin{enumerate}
		\item\label{thm:a1} \rew{\eqref{eqn:BQP} has no feasible solution
			with value} greater than $\rho$; and
		\item\label{thm:a2} for any $\x$ in the set $\Delta^c$, we have $h(\x) > \rho$.
	\end{enumerate}
        If  $f^* < \infty$
	then the optimal values of the constrained and the unconstrained problem coincide,
	i.e., $h^* = f^*$.
	Moreover, this instance is infeasible if and only if $h^* > \rho$.
\end{theorem}

The choice of parameters used in~\cite{GuWi:19} is
\begin{align*}
  \rho  &= \tilde{u}\\
  \sigma &= \tilde{u} - \tilde{\ell} + \epsilon,
\end{align*}
where
\begin{align*}
	\tilde{\ell} &= \min \left\{c^\top \x + \prods{F}{X} \mid \diag(X)=\e, \ X-\x\x^\top \succcurlyeq 0,  \mathcal{A}_{3}(X)\le \e ,  \mathcal{A}_{5}(X)\le \e   \right\} \\
	\tilde{u} &= \max \left\{c^\top \x + \prods{F}{X} \mid \diag(X)=\e, \ X-\x\x^\top \succcurlyeq 0, \
	[b, \unaryminus A ] \cdot
	\begin{bmatrix}
		1 & \x^\top \\
		\x & X
	\end{bmatrix} = \mathbf{0}\right\}.
\end{align*}
In this way, $\rho$ and $\sigma$ fulfill the assumptions of Theorem~\ref{thm:exp},
but are kept ``sufficiently'' small in order to avoid numerical difficulties when computing the maximum cut.
\rew{For a more detailled study on the choice of the parameters $\rho$ and $\sigma$ we refer to~\cite{GuWi:19}.}

In~\cite{GuWi:19}, further enhancements of the choice of the penalty parameter are discussed,
e.g., an update is made as soon as a feasible solution of \eqref{eqn:BQP} is found,
or an early stopping condition is added when infeasibility of \eqref{eqn:BQP} is detected.

\section{BiqBin solver}
\label{subsec:smc}

\subsection{(Sequential) branch-and-bound algorithm}
The BiqBin solver is a Max-Cut based solver for \eqref{eqn:BQP} instances that solves
their reformulation to \eqref{eqn:MC} instances using a B\&B algorithm.

The main ingredients of BiqBin are:
\begin{enumerate}
	\item the procedure for the exact penalty reformulation of \eqref{eqn:BQP} instances into \eqref{eqn:MC} instances;
	\item the bounding procedure, which provides for each instance of a problem (also for smaller subproblems obtained via branching) an upper bound on the optimum value;
	\item the branching procedure, which splits the current problem into more problems of smaller dimensions by fixing some variables;
	\item a heuristic for generating feasible solutions providing a lower bound.
\end{enumerate}

The overall performance of the algorithm
is determined by the quality of the lower and upper bounds, computed in each B\&B node,
by the computational efficiency of computing these bounds,
and by the strategy of exploring the B\&B tree.
The main reason why we cannot solve large instances on personal computers is that the B\&B tree grows too big,
i.e., the pruning is too slow and a single processor is not capable to explore all the generated nodes in the tree.

\paragraph{Exact penalty reformulation.} The procedure which reformulates every instance of \eqref{eqn:BQP} into an instance
of \eqref{eqn:MC} is  described in Section \ref{sec:exact_penalty}.

\paragraph{Bounding procedure.}
The starting point of the algorithm is the strengthened SDP relaxation \eqref{hyper_SDP}.
\rew{Since the number of inequalities is too large to solve this SDP by standard solvers, we will use a bundle method to find an approximate solution to the partial Lagrangian dual.}
By dualizing only the inequality constraints, we obtain the nonsmooth convex  partial dual function
$$
	f(\gamma) = \max_{\substack{\textrm{diag}(X) \ = \ \e,\\ X \succcurlyeq 0}} \mathcal{L}(X,\gamma)
	= \e^\top \vecgamma + \max_{\substack{\textrm{diag}(X) \ = \ \e,\\ X \succcurlyeq 0}}\langle L - \mathcal{A}_{\prob{HYP}}^\top (\gamma), X \rangle,
$$
where $\gamma$ are the nonnegative dual variables associated with the constraints $\mathcal{A}_{\prob{HYP}}(X)\le \e$.
Evaluating the dual function $f(\gamma)$ and computing the subgradient amounts to solving an SDP of the form \eqref{basic_SDP},
which can be efficiently computed using an interior-point method tailored for this problem.
It provides us with the matching pair $(X_{\gamma}, \gamma)$ such that $f(\gamma) = \mathcal{L}(X_{\gamma},\gamma)$.
Moreover, the subgradient of $f$ at $\gamma$ is given by $\partial f(\gamma) = \e - \mathcal{A}_{\prob{HYP}}(X_{\gamma}).$
For obtaining an approximate minimizer of the dual problem
\begin{align}\label{dual_problem}
	\begin{split}
	\min \quad& f(\gamma) \\
	\st \quad& \gamma \ge 0,
	\end{split}
\end{align}
we use the bundle method~\cite{Ki:89, ReRiWi:10}.
Let the current iterate be $\hat{\gamma}$.
Suppose we have evaluated $f$ at $k \ge 1$ points $\gamma_1,\ldots,\gamma_k$ with matching pairs $X_1,\ldots,X_k$ and subgradients
$\e - \mathcal{A}_{\prob{HYP}}(X_i)$ for $i\in \{1,\dots,k\}$.
The bundle method combines the following two ideas:
\begin{enumerate}
	\item the function $f(\gamma)$ is approximated by 
	\begin{align*}
		f_{\mathrm{appr}}(\gamma) &= \max\{\mathcal{L}(X,\gamma) \mid X \in \conv(X_1, \ldots, X_k)\}\\[0.5em]
		&= \max_{\lambda \ge 0, \ \e^\top \lambda = 1} \e^\top \gamma + \langle L-\mathcal{A}_{\prob{HYP}}^\top (\gamma),\mathcal{X}\lambda\rangle,
	\end{align*}
	where the bundle of matrices $\mathcal{X} = (X_1,\ldots,X_k)$ is used to construct a minorant $f_{\mathrm{appr}}$  and $\mathcal{X}\lambda=\sum_i \lambda_i X_i$;
	\item the proximal point idea, which penalizes the displacement from the current best point $\hat{\gamma}$
		with a quadratic regularization term proportional to $\Vert \gamma - \hat{\gamma}\Vert^2$.
\end{enumerate}
In summary, the bundle method finds a new trial point by minimizing, for some prescribed parameter $t > 0$, the function
$$
	f_{\mathrm{appr}}(\gamma) - \frac{1}{2t}\Vert \gamma - \hat{\gamma} \Vert^2
$$
over the nonnegative orthant.
As the method terminates, we obtain the approximate minimizer of the dual function,
as well as the convex weights $\lambda\ge 0,\sum_i \lambda_i=1$, determining the matrix $\mathcal{X}\lambda$.

In order to obtain a tight upper bound for an instance of \eqref{eqn:MC},
we use the cutting-plane approach,
where multiple $k$-gonal inequalities are added and purged in the course of running the bundle algorithm.
First, the optimum solution of the basic semidefinite relaxation \eqref{basic_SDP} is computed using an interior-point method followed by separating a set of triangle inequalities.
After doing a few bundle iterations with this set of constraints,
we purge all inactive constraints.
Next, invoking the separation routine described in Section~\ref{sec:sdp},
new violated $k$-gonal inequalities are added.
The problem with the new set of constraints is solved and the process is iterated as long as there is a significant decrease of the upper bound.

\rew{Since the primal solution matrix $X$ is not available,}
we use the same idea as in BiqMac to purge some inequality constraints. 
We look at the values of the corresponding dual multipliers in the vector $\gamma$.
If the value of some dual multiplier is close to zero, this indicates that the corresponding constraint is not active and we remove it.

After each iteration of computing an upper bound and separating triangle inequalities,
higher order $k$-gonal inequalities are added in order to further decrease the bound.
We monitor the maximum violation of triangle inequalities
$r_{\text{tri}}=\max(\mathcal{A}_3(\mathcal{X} \lambda ) - \e)$
and as soon as the number is sufficiently small,
the heuristic from Section~\ref{sec:sdp} is used to add some strongly violated pentagonal inequalities to the relaxation.
Similarly, as the maximum violation of pentagonal inequalities $r_{\text{pent}}$ drops below some threshold,
new heptagonal inequalities are separated and added to the relaxation.
In our numerical tests, we have used the thresholds $r_{\text{tri}} < 0.2$ and $r_{\text{pent}}<0.4$.

In order to improve the performance of BiqBin,
we can stop the bounding routine when we detect that we will not be able to prune the current node in the B\&B tree.
We again borrow an idea from BiqMac.
After some cutting plane iterations, we make a linear (and hence optimistic) forecast to decide whether it is worth doing more iterations.
If the gap cannot be closed, we terminate the bounding routine, branch the current node, and start evaluating new subproblems.
This is especially important in the parallel solver, since its efficiency depends on how quickly idle workers receive subproblems.

Furthermore, at the beginning of the algorithm, the number of bundle iterations should be small.
We start with three iterations.
Then this number is increased after each separation of new cutting planes, until a limit is reached.
In the BiqBin solver, this limit value is set to~15.
This is motivated by the fact, that in the beginning it will take a while until we have identified the right cutting planes
and we do not waste time by trying hard to decrease the bound if the current set of hypermetric inequalities does not allow much progress.

\paragraph{Branching strategies.}
BiqBin uses two branching strategies which are based on the bundle of matrices obtained from the bounding routine.
Once the bundle method terminates, the last column $\x$ of the matrix $\mathcal{X}\lambda$ is extracted.
Due to the diagonal and positive semidefiniteness constraints on the feasible matrices of \eqref{hyper_SDP},
all the entries of $\x$ lie in the interval $[\unaryminus1,1]$.
Similarly to BiqCrunch, the decision on which variable $x_i$ the branching should be performed is based on the following two strategies:

\begin{enumerate}
	\item difficult first: we branch on the vertex $i$ for which the variable $x_i$ is closest to $0$;
	\item easy first: we branch on the vertex $i$ for which the variable $x_i$ is furthest from $0$.
\end{enumerate}
In the $0/1$ formulation, these rules are usually referred to as most-fractional and least-fractional rules, respectively.
The difficult first rule is set as the default strategy for the BiqBin solver.

When branching, two new subproblems, in which the chosen branching variable is fixed accordingly,
are created and the corresponding nodes are added in the B\&B tree, i.e., to the priority queue of unexplored problems.
Priority is based on the upper bound obtained from the bundle method.
When selecting the next subproblem, a node with the worst upper bound is evaluated first.

\paragraph{Rounding heuristic.}
For generating high quality feasible solutions of the Max-Cut problem,
we apply the Goemans-Williamson rounding hyperplane technique~\cite{GoWi:95}.
Let $X$ be an optimal solution of some SDP relaxation of Max-Cut.
By computing the Cholesky factorization $X = V^\top V$ with column vectors $v_i$ of $V$ and selecting some random vector $r$,
the cut $(S,V\backslash S)$ is obtained by setting
$
	S = \{i \mid v_i^\top r \ge 0 \}.
$
Since in our case we are working with the partial Lagrangian,
the information about the primal $X$ is not available.
Instead, we use a convex combination of bundle matrices $\mathcal{X}\lambda$ as the input.
The cut vector $x$ obtained from this heuristic is then further improved by flipping the vertices
and using a convex combination of $\mathcal{X}\lambda$ and the cut matrix $xx^\top$.
To summarize, for generating good cuts, we use the following iterative scheme:
\begin{enumerate}
	\item use the Goemans-Williamson rounding hyperplane technique to generate cut vector $x$ from  $\mathcal{X}\lambda$;
	\item the cut $x$ is locally improved by checking all possible moves of a single vertex to the opposite partition block;
	\item by using a convex combination of $\mathcal{X}\lambda$ and $xx^\top$,
		we bring the rounding matrix towards a good cut.
		With this new matrix, go to step~(1)
		and repeat as long as one finds better cuts.
\end{enumerate}
Interestingly, in most of our numerical experiments, this heuristic finds the optimum solution already in the root node.

\paragraph{Strategy for faster enumeration of the B\&B tree.}
As described in Section~\ref{sec:sdp},
simulated annealing is used to heuristically separate $k$-gonal inequalities.
Adding these $k$-gonal inequalities to the model in each B\&B node is not necessary,
especially if one can not prune the node.
In that case all the work done and time are wasted,
since for bigger graphs we usually need to reach a certain depth in the B\&B tree in order to prune the nodes,
even if the bounding routine produces tight bounds.
It is beneficial to check
(before including cutting planes when processing the node)
whether there is hope to prune the node or whether it is better to branch and produce smaller subproblems.
We propose the following strategy.

In the root node we compute:
(i) the bound $\OPT_{\prob{SDP}}$ of the basic SDP relaxation \eqref{basic_SDP}, which is not strong but is quick to compute; and
(ii) the bound $\OPT_{\prob{HYP}}$ by iteratively including violated triangle, pentagonal, and heptagonal inequalities.
Let
\label{page:diff}
\begin{equation}
	\label{eq:diff}
	\diff = \OPT_{\prob{SDP}} - \OPT_{\prob{HYP}},
\end{equation}
and let $\LB$ denote the current lower bound.
Then, at all other nodes, we first compute only the basic SDP bound $\OPT_{\prob{SDP}}$.
If the condition
\begin{equation}
	\label{diff}
	\OPT_{\prob{SDP}} \le \LB + \diff + 1
\end{equation}
is satisfied,
this means, we are already close to the lower bound, and so we add cutting planes to compute the
tighter bound $\OPT_{\prob{HYP}}$
in order to increase the probability of pruning the node.
With this idea, we can efficiently traverse the B\&B tree and only invest time into the bounding routine when it is needed.
Numerical results show that overall this strategy produces more B\&B nodes than necessary, but the performance of the algorithm improves;
in particular, it has a positive impact on the parallel version described in Section~\ref{para_bab}.

\paragraph{Numerical illustration of the bounding procedure.}	
When using a B\&B algorithm, one has two options regarding the quality of the upper bound.
On one hand, a strong upper bound can be computed by iteratively adding and purging multiple cutting planes.
In this way, more work is done within a node but overall this approach produces fewer B\&B nodes.
On the other hand, one can efficiently compute slightly weaker upper bounds,
hence the whole tree grows larger but if the time spent for evaluating each node is small it can be traversed faster. In BiqBin we use the first approach.

\begin{figure}[ht!]
	\centering
	\includegraphics[scale=0.6]{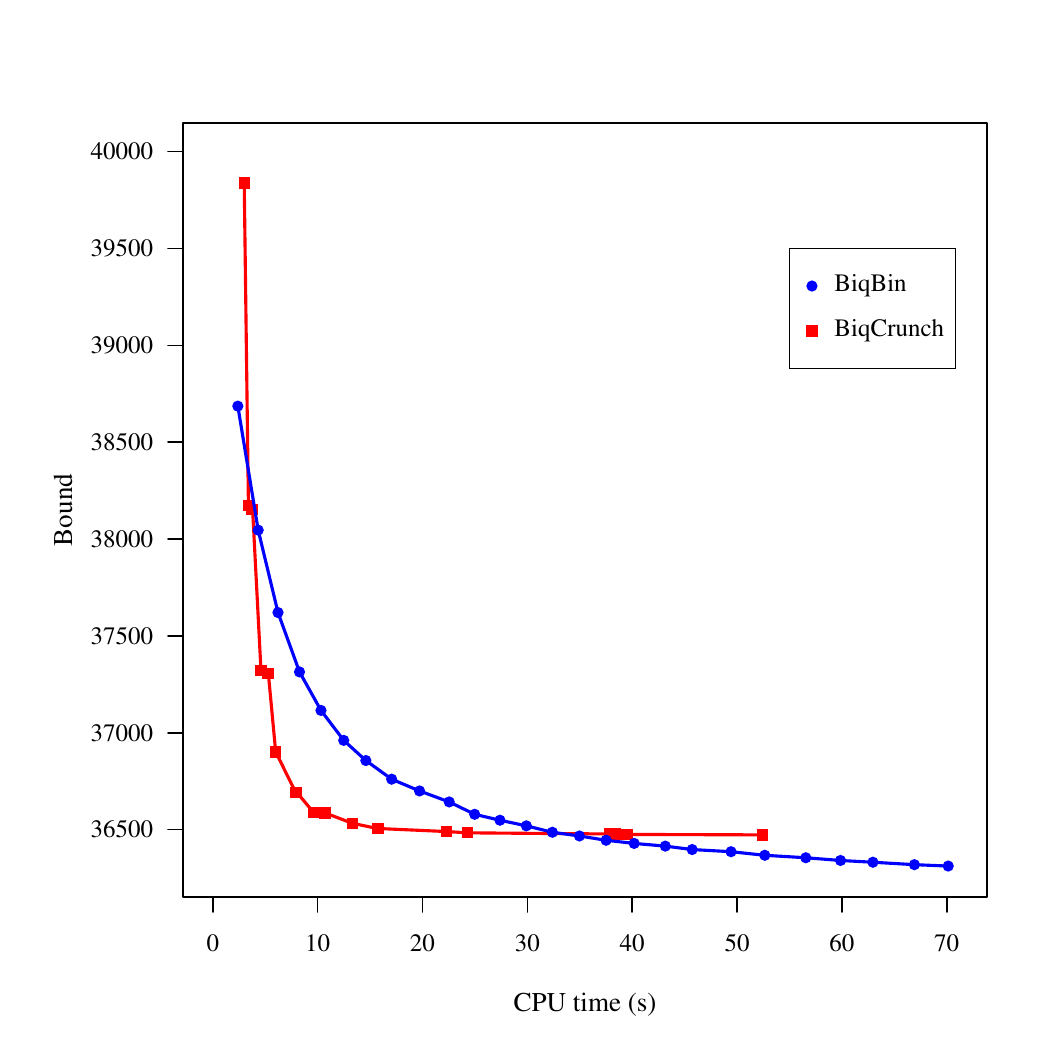}
	\caption{Comparison of the bounding procedure of BiqCrunch and BiqBin in the root node on the Beasley bqp250.8 problem.
		The upper bound computed with BiqCrunch is approx.~$36\,472$, whereas with BiqBin we obtain approx.~$36\,287$. }
	\label{fig:decrease_bound}
\end{figure}

We take the Beasley bqp250.8 problem from the BiqMac library and plot the convergence curve
for the bounding routines in the solvers BiqBin and BiqCrunch in the root nodes of the B\&B trees.
Figure~\ref{fig:decrease_bound} depicts the decrease of the dual function values in the course of a bound computation.
We note that by adding higher order $k$-gonal inequalities,
our bounding routine attains a tighter bound compared to BiqCrunch.
Consequently, BiqBin creates a smaller B\&B tree, which  consists of 81~nodes, while BiqCrunch terminates after traversing 325~nodes.

\subsection{Parallelization of branch-and-bound}\label{para_bab}
In this \rew{section}, we describe how the algorithmic ingredients from sequential B\&B algorithm are combined
into a parallel solver which utilizes distributive memory parallelism.
This is done in a similar fashion as in the parallel B\&B solver for the stable set problem~\cite{Hrga_et_al_2020}, introduced by some of the authors of this paper.

The load coordinator--worker paradigm with distributive work pools is applied, in which the rank~0 process becomes the master process
carefully managing the status of each worker (idle or busy),
while different workers concurrently explore branches of the B\&B tree.
Each worker has its own local queue of subproblems and the work is shared when one of them becomes idle.
The master node knows the status of each worker and acts as a load coordinator receiving messages and, based on their content,
replying in an appropriate manner.

At the beginning of the algorithm, the load coordinator reads and broadcasts the original graph to the workers and initializes the solution.
It is important that every process has the knowledge of the original graph,
since construction of subproblems via branching and encoding the MPI messages is done based on this information.
All the data about the B\&B nodes is encoded as an MPI structure,
which is used in communication between different workers in order to efficiently exchange and construct the subproblems.

Next, the master process evaluates the root node and distributes the best lower bound.
After the bounding step, two new subproblems are generated,
which are sent to the first two idle processes.
Afterwards, its job is restricted to monitoring the status of the workers,
counting the number of B\&B nodes,
and distributing the best solution found so far.
In order to solve multiple SDP relaxations in parallel,
the workers need to send and receive subproblems and keep load balance.
In this way, the work is shared and the whole dynamic B\&B tree is enumerated faster.

After the initialization phase, the master process waits for three types of messages sent by the workers.
Firstly, if the worker's local queue is empty, the message is sent informing the master that the process is idle and can receive further work.

Secondly, the master process receives messages regarding the load balance and sharing of work.
Throughout the algorithm, the master node manages a bool array containing the statuses of the workers.
They specify whether the process is active or idle.
The load coordinator uses this information to reply to the requests received from the workers during the branching step,
in which the process wants to share one of the newly generated subproblems with some other free worker.
The load coordinator sends a message specifying which workers are free and the value of best lower bound.

And thirdly, during the execution of the algorithm, the working processes compute multiple candidates for optimum solutions.
The master node is keeping track of the currently best value and the corresponding solution.
When a new solution is received, the value is compared, updated if necessary, and distributed back during the communication phase.

After a worker computes lower and upper bounds,
it compares these values to see whether this branch of enumeration tree can be safely pruned or further branching is needed
and construction of new subproblems takes place.
In the latter case, a request message is sent to the load coordinator,
asking for idle processes to share one of the newly generated subproblems or subproblems left in the queue from previous branching processes.
If no idle worker is available, the generated subproblems are placed in the worker's queue and the work continues locally.
Otherwise, subproblems are encoded and sent to available idle workers.
This is also where the exchange of the best lower bound happens.

When all the workers become idle,
the master process sends a message to finish and the algorithm terminates.
The algorithms for the load coordinator and the workers are summarized in Algorithms~\ref{alg:master} and~\ref{alg:worker}.
\begin{algorithm}[ht!]
	\caption{Algorithm for load coordinator}
	\label{alg:master}
          \KwData{Graph $G=(V,E)$.}
           read the input file and distribute the graph $G$ to the workers\;
	   evaluate the root node\;
           \If{not pruned}{
             broadcast the current best lower bound and variable $\diff$\;
             set status of every worker to \emph{idle}\;
             branch root node, generate subproblems and send them to first two idle workers; set their status to \emph{busy}\;
             \Do{busy workers exist}{
             receive message from worker\;
             \uIf{idle}{
               mark worker as \emph{idle}\;}
             \uElseIf{new\_value}{
                 worker found better lower bound $\LB$\;
                 compare, update and return the new value\;}
             \uElseIf{send\_workers}{
                   send ranks of the idle workers to the sender\;}
                 }
	     send \emph{finish} message to the workers\;}
\end{algorithm}

\begin{algorithm}[ht!]
  \caption{Algorithm for worker process}
  \label{alg:worker}
  receive graph $G$, best lower bound $\LB$ and $\diff$ from the load coordinator\;
		 \Do{not finish}{
		Receive \emph{finish} message or \emph{subproblem} from the load coordinator or other worker and put the problem into the local queue\;
		\If{subproblem received}{
		\While{queue is not empty}{
			select the next subproblem and compute lower bound $\mathrm{lb}$ and upper bound $\mathrm{ub}$\;
			\If{$\mathrm{lb} > \LB $}{
				send \emph{new\_value} message and $\mathrm{lb}$ to the load coordinator and receive updated $\LB$\;}
			\If{$\LB < \mathrm{ub}$}{
				Generate new subproblems and place them in local queue\;
				Send request message \emph{send\_workers} and receive ranks of idle workers from load coordinator\;
				Encode subproblems and send them to available processes\;
				}
			}
		Send \emph{idle} message to the load coordinator and wait for new tasks.
		}
		}	
	\end{algorithm}

Lastly, we explain how the parallel version benefits from the strategy using the variable ``$\diff$''
(see equation~\eqref{eq:diff} on page~\pageref{page:diff}).
If the number of available workers is large,
we need to reach a certain depth in the B\&B tree in order for the processes to receive the work.
Until this happens in the algorithm, the workers are idle.
To fully exploit all HPC resources available, we need a strategy for the worker processes to start evaluating the nodes as soon as possible.
This is where the idea from Section~\ref{subsec:smc} helps.
After the load coordinator evaluates the root node,
the best lower bound found and the variable $\diff$ are distributed to all the workers.
When the first two idle processes evaluate the generated subproblems,
typically the value of the basic SDP relaxation is such that condition~\eqref{diff} is not satisfied.
Hence, on the first few levels of the B\&B tree, the workers compute only the basic SDP bound.
Then branching of new subproblems takes place and idle workers quickly receive the generated subproblems.
After the bounds are such that condition~\eqref{diff} is valid,
hypermetric inequalities are added to compute tighter bounds $\OPT_{\prob{HYP}}$ approximately.
This implies that more nodes are pruned, meaning that the size of the B\&B tree decreases, and thus the algorithm terminates faster.

\section{BiqBin web application}

 The BiqBin solver is available as a web application.
 Its main purpose is to enable (registered) users to test the solver on their own \rew{problem} data.

 The web application contains the main information about the BiqBin project
 and benchmark results, which are publicly available.
 Upon registration, a user can access the core functionalities of the application.
 In short, the users upload their problem data files and decide which algorithm/solver they want to use to solve the problem.
 The application then sends these data to an HPC,
 where the selected algorithm is executed against the given data,
 and, after some time interval, it returns the solution (optimum or approximate)
 back to the web application,
 which finally notifies the users by an email that the result is ready.

 Let us describe the web application in more details.
 The initial step for a user is to prepare the problem data in an appropriate format,
 which is described in the description of each \rew{function\footnote{\url{http://biqbin.eu/Home/Features}}.}
 Then, the user creates a new \textit{instance}, which is the entity determining the problem data.
 The user determines the instance's name, data (by uploading the file with all the data),
 and description. In this way, the user is able to reuse the data later, e.g.,
 for testing the performance of other solvers on the same instance.

 Next, the user creates a \textit{task}, which is comprised of an instance and a solver (function)
 the user wants to run with the instance's data.
 The user also determines the task's name,
 and specifies if the result and the instance's data can be listed publicly for benchmarking purposes.
 When all the properties are specified, the task can either be started automatically,
 or the user can decide to start it later.

 \textit{Starting a task} means that the system sends a demand with the execution data to a specified HPC.
 The administrators can configure which HPC should be used for each of the available functions
 or even for a user.
 For each sent task, a job is started on a specified number of cores.
 Execution of a solver is being regularly monitored by saving intermediate results per every given time interval.
 After a specified maximum execution time, or after finding an optimum solution,
 the solver's execution is stopped, and the result files are sent back to the web application together with the task's metadata.
 The result is saved in the database and the user is notified by an email that the job is completed.
 If the task has been marked as publicly available, then a benchmark record is also created and immediately visible on Benchmarks site.

 The system is designed in such a way that additional solvers or algorithms can be easily added\rew{.}
 The number of cores used for a single task is fixed, so that one user does not use all of the available infrastructure.
 However, the users, who have special roles assigned, can specify the number of cores themselves.
 In this way, the architecture of the system is highly extendible in all main directions.
 Scalability depends on the HPC facilities, but this is also considered in the web application,
 since tasks for every solver can connect to a selected HPC.
 \begin{figure*}[htp]
 	$$
 	\includegraphics[scale=0.7]{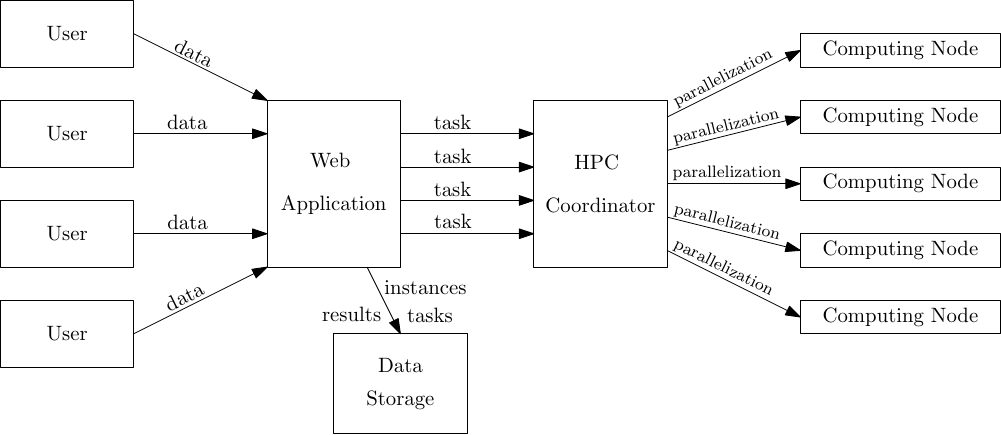}
 	$$
 	\caption{A workflow diagram from users' inputs to algorithms' execution on HPC} 
 	\label{fig:schema}
 \end{figure*}

\section{Numerical results}

In the numerical part of the paper,
we benchmark the BiqBin solver with BiqCrunch~\cite{krislock2017biqcrunch}, Gurobi~\cite{gurobi}, and SCIP~\cite{achterberg2009scip}
on four families of test instances of \eqref{eqn:BQP}, which we describe in the following paragraphs.
We decided for these solvers because they are well-known to be among the best solvers for these type of problems and are freely available for research purposes.
\rew{All the solvers received as input the original problem~\eqref{eqn:BQP},
hence only BiqBin applies the reformulation into a max-cut instance.
This means that Gurobi and SCIP use all the MIP techniques that are available.}
We ran Gurobi and SCIP using the default parameters.
Performance of the BiqCrunch solver depends on the branching rule used.
Thus we decided to use both branching rules,
i.e., most-fractional and least-fractional rule, and we report the results for the faster run.
For other paramers of BiqCrunch, we used their default values.
We also planned to include CPLEX in this study but we were not able to obtain the academic licence for the HPC system that we were using.

All computations were performed on the HPC system at University of Ljubljana, Faculty of mechanical engineering.
We used the E5-2680 V3 (1008 hyper-cores) DP cluster, with IB QDR interconnection,
164 TB of LUSTRE storage, 4.6 TB RAM, supplemented by GPU accelerators.
As a Message Passing Interface library we used Open MPI and the code is compiled against OpenBLAS and LAPACK.
The BiqBin solver was run directly on the cluster, so the computational times might be slightly different
compared to the times obtained by the BiqBin on-line solver.

\subsection{Four special families of \eqref{eqn:BQP} instances used for numerical tests}

\paragraph{Max-Cut instances}
The first family of \eqref{eqn:BQP} benchmark instances consists of Max-Cut instances.
\rew{We selected benchmark instances from the BiqMac library that have been considered as hard by other solvers} (see~\cite{biqmaclib} for more details).
Furthermore, we used rudy, a graph generator written by Giovanni Rinaldi~\cite{Ri},
to construct new (hard) Max-Cut instances with 180 nodes, similar to the ones in the BiqMac library.

\paragraph{Unconstrained BQP instances}
The second family of \eqref{eqn:BQP} instances consists of two sets of unconstrained BQP.
The instances of the first set, due to Beasley~\cite{Be:98}, have density 0.1 with values uniformly chosen from $[\unaryminus 100,100]$.
The second set of unconstrained BQP instances were generated by Billionet and Elloumi~\cite{BiEl:07} according to~\cite{PaRo:90}.
They have size $n \in \{100,120,150,200\}$ and different densities.
The diagonal coefficients are in the range $[ \unaryminus 100,100]$, while the off-diagonal ones are in from $[ \unaryminus 50,50 ]$.
Similar instances with size~250 and density~0.1 were taken from the BiqMac library.
For more details on these instances we refer to~\cite{biqmaclib}.

\paragraph{Densest $k$-subgraph instances}
The third family of benchmark problems consists of instances of the \eqref{eqn:DkS},
described in Section~\ref{sec:mot}.
It is a binary quadratic problem with one linear constraint, also called the cardinality constraint.
\rew{These problems, also known as cardinality Boolean quadratic programming problems, include a wide range of applications in telecommunications or chemistry, see~\cite{lima-grossmann2017} for further details.}
We solve the benchmark instances that can be found at~\cite{kcluslib}.
They have different sizes $n \in \{120,140,160\}$, densities $d \in \{ 0.25, 0.50, 0.75\}$, and values for the parameter $k \in \{ n/4, n/2, 3n/4\}$.
For testing the parallel version of BiqBin, we also created similar instances with sizes $n \in \{180,200\}$.
They can be found at the BiqBin web page\footnote{\url{http://biqbin.eu/Home/BenchmarkInstances}}.

\paragraph{Randomly generated instances of \eqref{eqn:BQP}}
The fourth family of instances consists of randomly generated instances with a varying number of linear constraints.
These instances have size $n = 100$ and up to 15~constraints, their description can be found in~\cite{GuWi:19}
and they are available at~\cite{Gu:19}.

\subsection{Comparison of sequential algorithms}

In this section, we compare the sequential version of BiqBin with BiqCrunch~\cite{krislock2017biqcrunch}, GUROBI~\cite{gurobi}, and SCIP~\cite{achterberg2009scip}
solvers on the four families of problems introduced in the previous section.
We note that for the first three families of benchmark problems GUROBI and SCIP did not solve any instance within 3~hours,
hence we report results for these problems only for BiqBin and BiqCrunch.

We present the results of the comparisons of these two solvers on the Max-Cut instances in Tables~\ref{tab:BiqBin_MC}--\ref{tab:BiqCrunch_MC},
the results for other unconstrained BQPs are in Tables~\ref{tab:BQP_biqBin}--\ref{tab:BQP_BiqCrunch},
and the results for the densest $k$-subgraph problem are in Tables~\ref{tab:densestBiqBin}--\ref{tab:densestBiqCrunch}.

\begin{figure}
 \begin{subfigure}[b]{0.5\textwidth}
                \includegraphics[scale=0.4]{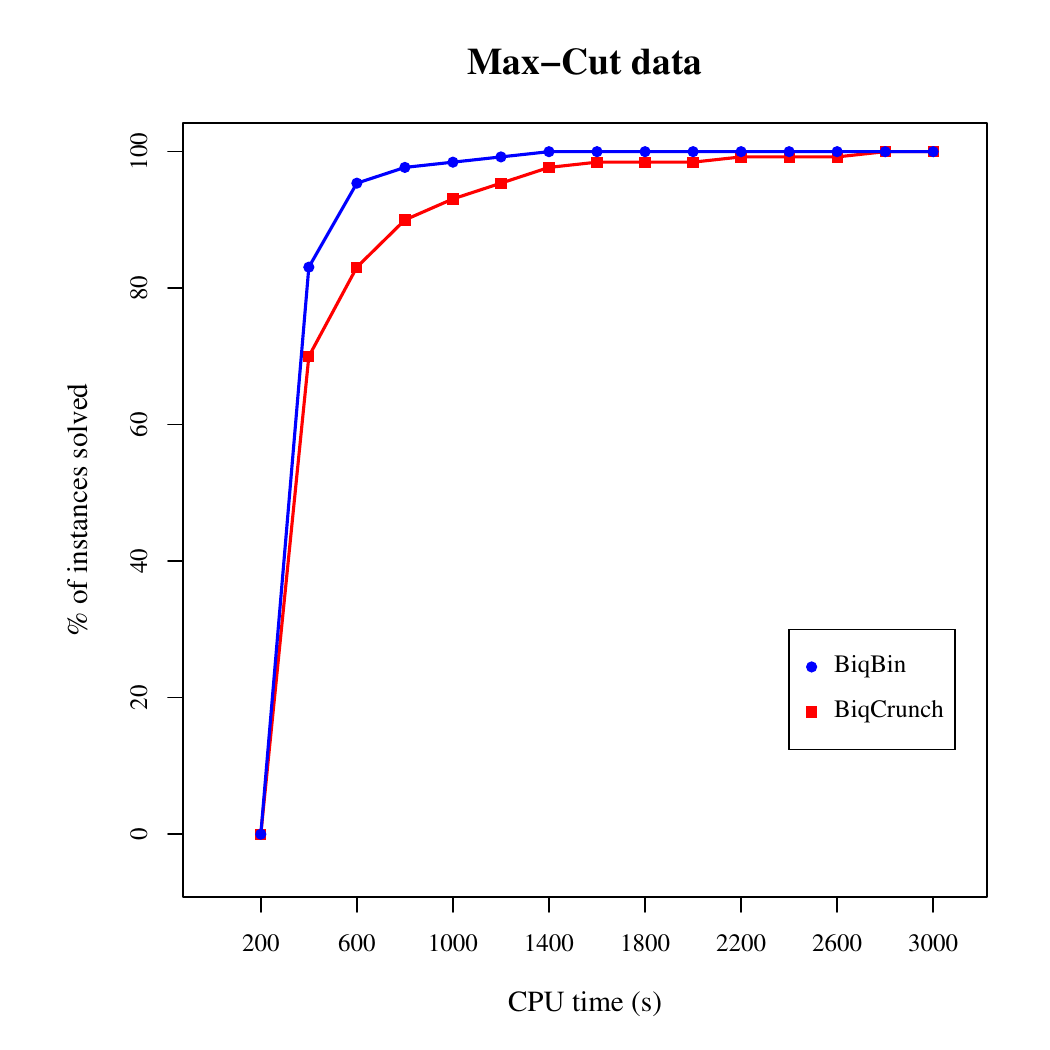}
                \caption{BiqBin vs.~BiqCrunch on the Max-Cut instances}
        \end{subfigure}%
  \begin{subfigure}[b]{0.5\textwidth}
                \includegraphics[scale=0.4]{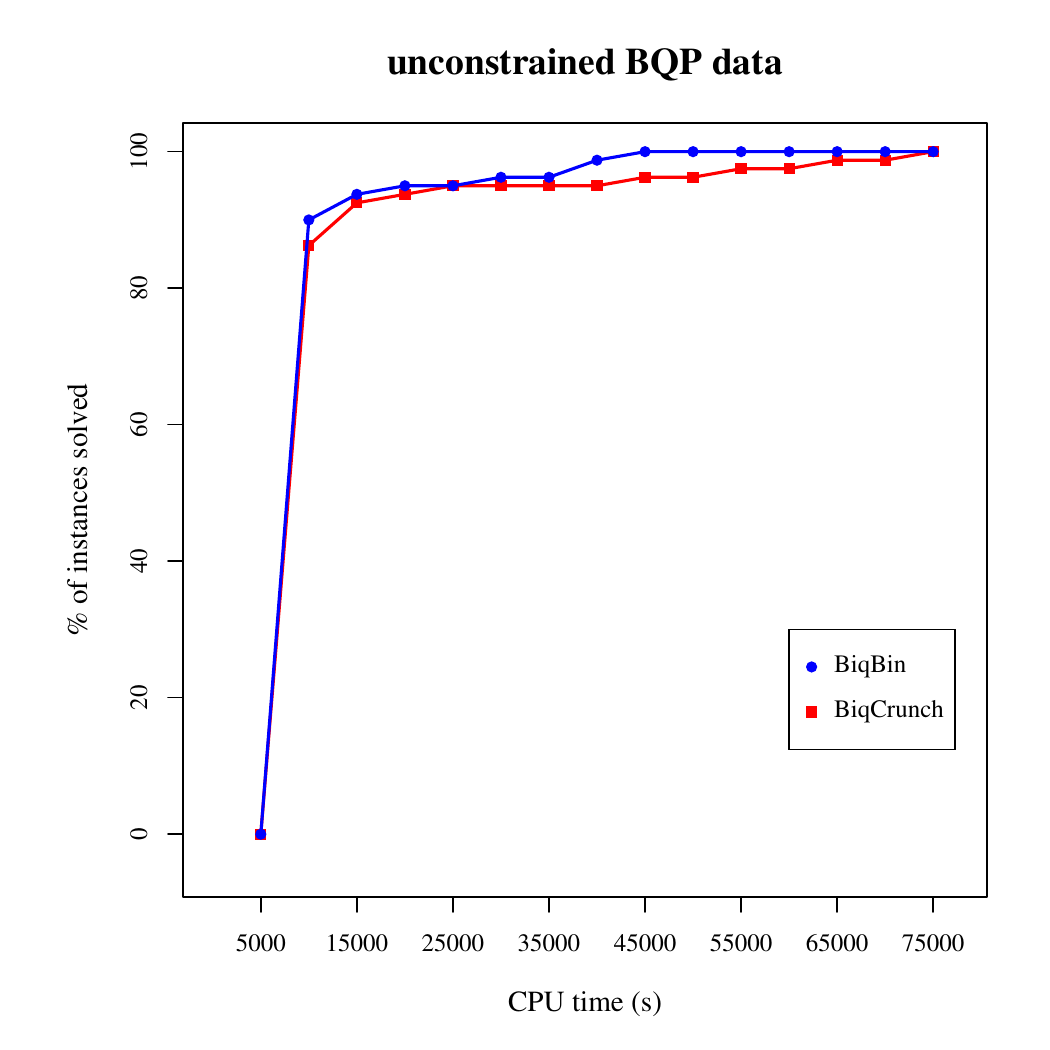}
                \caption{BiqBin vs.~BiqCrunch on the UBQP instances}
        \end{subfigure}\\
\begin{subfigure}[b]{0.5\textwidth}
                \includegraphics[scale=0.4]{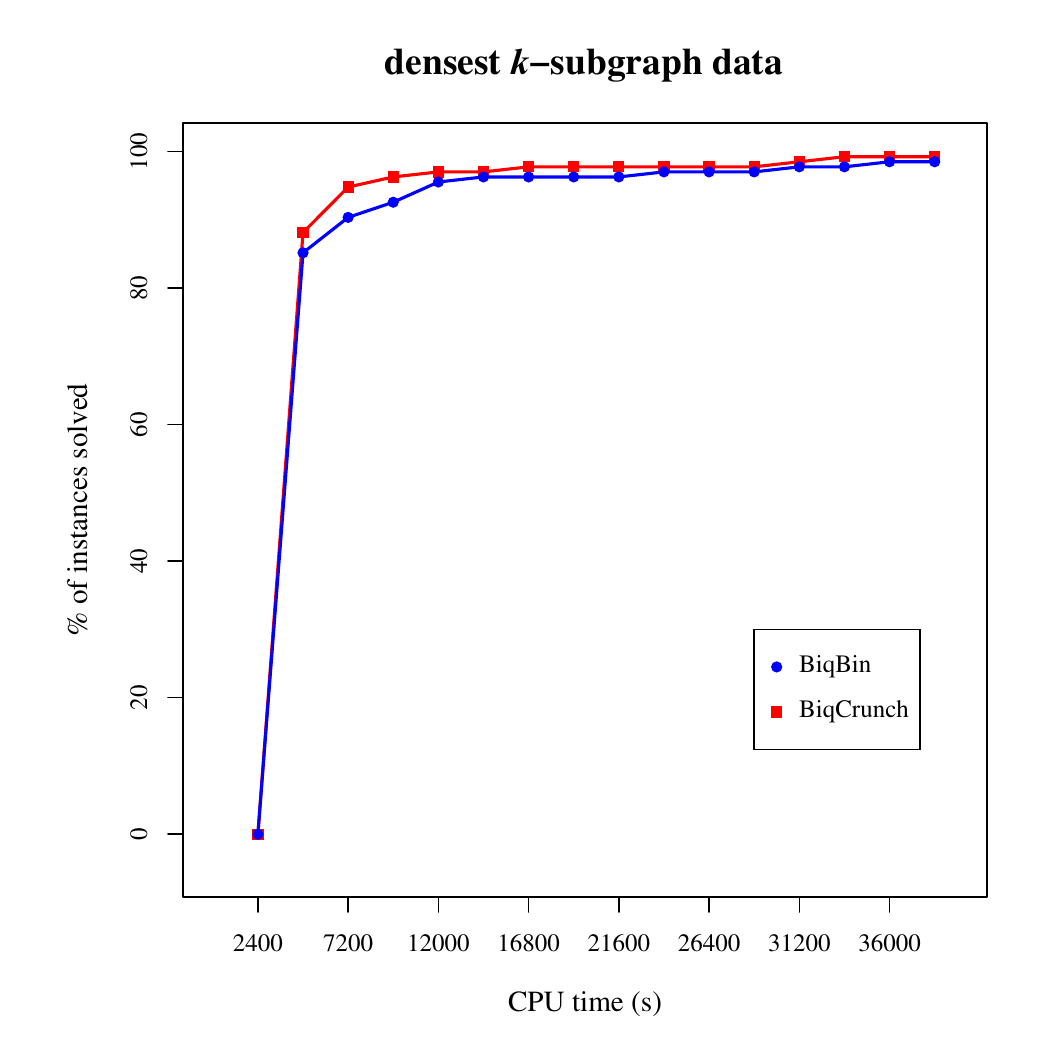}
                \caption{BiqBin vs.~BiqCrunch on the densest \\ $k$-subgraph  instances}
        \end{subfigure}%
  \begin{subfigure}[b]{0.5\textwidth}
                \includegraphics[scale=0.4]{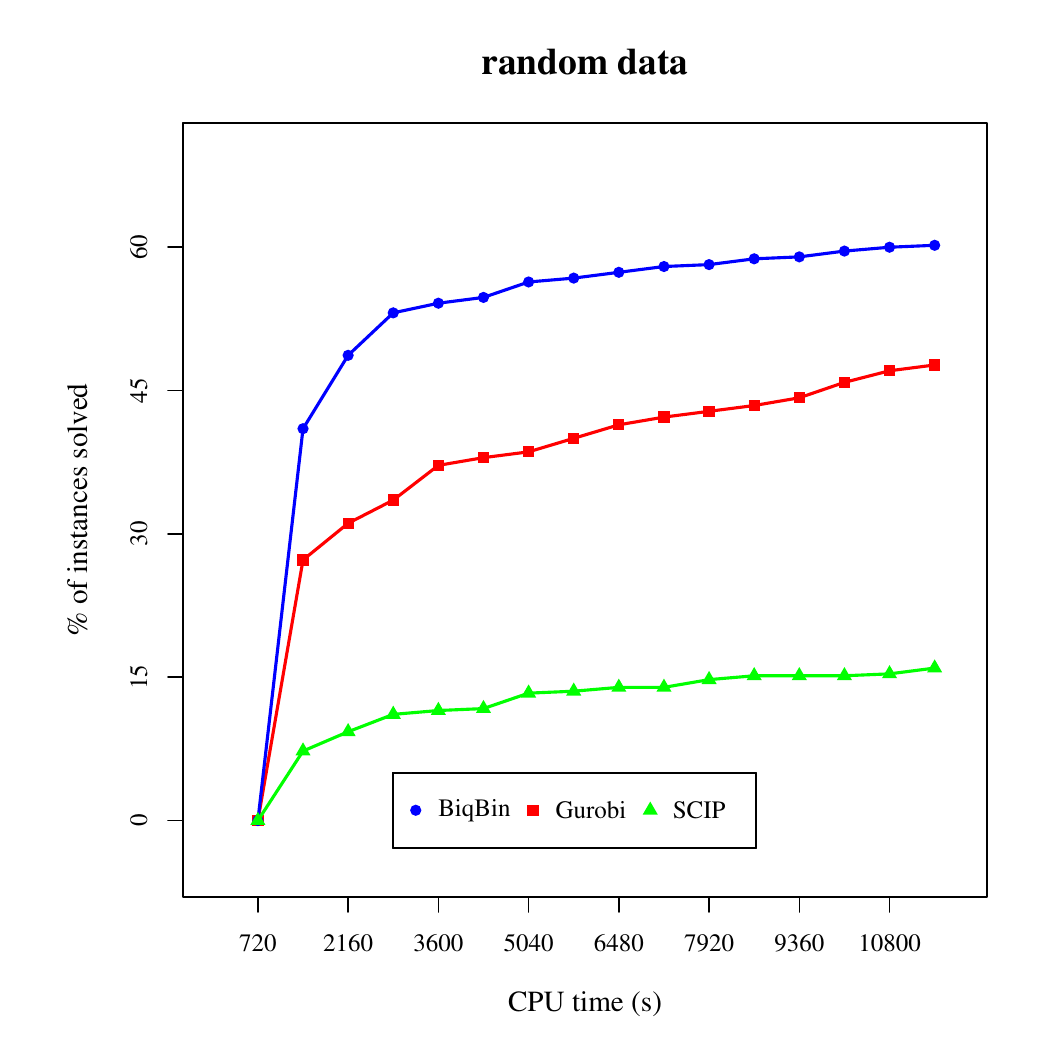}
                \caption{BiqBin vs.~Gurobi vs.~ SCIP on the randomly \\ generated instances of BQP}
        \end{subfigure}%
\caption{Performance profile for BiqBin and BiqCrunch on the Max-Cut, unconstrained Binary Quadratic Problem and densest $k$-subgraph data (Figures a, b, and c), and for solvers BiqBin, Gurobi and Scip on random data (Figure d)}\label{fig:per_prof}
\end{figure}

\rew{Before we list the tables, we summarize the results by the performance profile diagrams on Figure \ref{fig:per_prof}. These diagrams were created using the same raw data as we summarised in Tables  \ref{tab:BQP_biqBin}--\ref{tab:RND_biqBin}. There were large time spans for each solver and each data set, so the diagrams are quite coarse-grained. Nevertheless, we can see that on the Max-Cut, unconstrained BQP and randomly generated instances of BQP our solver BiqBin is outperforming the other solvers. On the instances of the densest $k$-subgraph problem BiqCrunch is performing slightly better then BiqBin.}

The first three columns in each of these tables contain the data about the instances:
for each group of instances, we used several instances,
e.g., for the Max-Cut problem and for the unconstrained BQP, we used 10 instances for every combination of size $n$ and density $d$,
and for the densest $k$-subgraph we used 5 instances.
The second set of columns in each table contains the minimum, the maximum, and the average number of nodes in the B\&B trees generated by BiqBin and BiqCrunch, respectively.
Similarly, the third set of columns reports the minimum, the maximum, and the average computing times \rew{(in seconds)} needed for the instances,
and the last set of columns contains the minimum, the maximum, and the average initial gaps (in \%).
\label{tab1:main_text}

For the Max-Cut problem, we can see in Tables~\ref{tab:BiqBin_MC}--\ref{tab:BiqCrunch_MC}
that the (sequential) BiqBin is always approximately 2~times faster than BiqCrunch and produces much fewer nodes in the B\&B tree.
One reason for this lies in the fact that the initial gap in the root node
(the difference between the upper and lower bound, divided by the lower bound) is smaller for BiqBin.

For the unconstrained BQP (Tables~\ref{tab:BQP_biqBin}--\ref{tab:BQP_BiqCrunch}),
we can see that both, BiqBin and BiqCrunch, perform very similar.
However, for BiqBin the B\&B trees are much smaller, so there is potential for improving the computing times by increasing the computational efficiency in each B\&B node.

For the densest $k$-subgraph problem (Tables \ref{tab:densestBiqBin}--\ref{tab:densestBiqCrunch}),
we can observe that BiqCrunch is in most cases slightly better (sometimes also much better) than BiqBin.
The reason for this is that BiqCrunch has a specific version adapted to solving this problem,
where the basic semidefinite relaxation is re-enforced with triangle and product constraints.

Turning to the randomly generated instances, GUROBI and SCIP performed well,
while the BiqCrunch solver showed very poor performance.
Thus we report in Table~\ref{tab:RND_biqBin} only the results for BiqBin, GUROBI and SCIP.
The randomly generated instances have different intervals from which the entries of $A$, $b$, and $F$ are selected,
thus for every combination, we have~15 instances with an increasing number of constraints, from~1 to~15.
Table~\ref{tab:RND_biqBin} contains in the first column the name of the family of instances,
the next three columns contain data about the numbers of instances (always 15),
the sizes of problem $n$ (always 100), and the densities of the data (the number of non-zero elements in $F$ divided by $n^2$).
For each of the three solvers, we report four columns containing the numbers of instances solved by the solver in the time limit of 3 hours,
and the minimum, the maximum, and the average computation times needed by particular solver on the instances solved within the time limit, respectively.

From Table~\ref{tab:RND_biqBin} it can be observed that{\rew,} overall, BiqBin is outperforming all the other solvers.
It solves for each family at least 2 out of the 15 instances,
and on average it solves 60\% of the instances, while GUROBI and SCIP solve 48\% and 16\% of the instances, respectively.
The computation times also demonstrate that BiqBin is on average much faster on the instances that it solves.

We can see that BiqBin has worst performance on instances starting with 100\_A\_0\_3\_ in the middle part of the table.
We do not have an unambiguous explanation for this, but the reason is very likely that the linear constraints enable generating good cutting planes, thus helping GUROBI and SCIP.
Obviously, the linear constraints for the instances  starting with 100\_A\_0\_3\_ yield more efficient cutting planes which enable these two solvers
to more efficiently prune the B\&B tree and therefore more often finish in a given time limit.
Some of the 100\_A\_0\_3 instances are also infeasible, while the instances having in the middle \_b\_0\_ are always feasible (they have zeros on the right hand, i.e., $\b = \reww{\mathbf{0}}$,
hence $\x = \reww{\mathbf{0}}$ is always a feasible solution). GUROBI and SCIP are solving problems in the \eqref{eqn:BQP} form, while BiqBin first reformulates the instances into the \eqref{eqn:MC} format.
This suggests to improve the performance of BiqBin by enhancing  an early  stopping condition for the case of  infeasibility using Theorem  \ref{thm:exp}.


\begin{landscape}
	\thispagestyle{empty}

\begin{table}[h!]\scriptsize
\centering
	\begin{tabular}{|l||rrr|rrr|rrr|rrr|}\hline
\shortstack{instance \\ group} 	 & 	 \# inst. 	 & 	 $n$ 	 & 	 density 	 & 	 \shortstack{B\&B \\ (min)} 	 & 	 \shortstack{B\&B \\ (max)} 	 & 	 \shortstack{B\&B \\ (avg)} 	 & 	 \shortstack{time \\ (min)} 	 & 	 \shortstack{time \\ (max)} 	 & 	 \shortstack{time \\ (avg)} 	 & 	 \shortstack{init. gap \\ (min)} 	 & 	 \shortstack{init. gap \\ (max)} 	 & 	 \shortstack{init. gap \\ (avg)} \\\hline\hline
g05\_100	 & 	10	 & 	100	 & 	0.50	 & 	11	 & 	697	 & 	191.8	 & 	18.7	 & 	951.2	 & 	\textbf{196.1}	 & 	0.3\% 	 & 	1.1\% 	 & 	0.7\%\\\hline
pm1d\_100	 & 	10	 & 	100	 & 	0.99	 & 	21	 & 	839	 & 	282.4	 & 	40.0	 & 	681.4	 & 	\textbf{231.1}	 & 	2.4\% 	 & 	8.1\% 	 & 	4.9\%\\\hline
pm1s\_100	 & 	10	 & 	100	 & 	0.10	 & 	1	 & 	15	 & 	5.8	 	 & 	0.8	 	 & 	22.5	 & 	\textbf{8.1}	 & 	0.7\% 	 & 	3.3\% 	 & 	1.5\%\\\hline
pw01\_100	 & 	10	 & 	100	 & 	0.10	 & 	1	 & 	13	 & 	4.2	 	 & 	1.5	 	 & 	31.4	 & 	\textbf{12.5}	 & 	0.0\% 	 & 	0.5\% 	 & 	0.1\%\\\hline
pw05\_100	 & 	10	 & 	100	 & 	0.50	 & 	13	 & 	289	 & 	112.2	 & 	53.5	 & 	503.9	 & 	\textbf{204.2}	 & 	0.1\% 	 & 	0.9\% 	 & 	0.6\%\\\hline
pw09\_100	 & 	10	 & 	100	 & 	0.90	 & 	39	 & 	201	 & 	114.4	 & 	67.9	 & 	349.4	 & 	\textbf{187.7}	 & 	0.2\% 	 & 	0.5\% 	 & 	0.4\%\\\hline
w01\_100	 & 	10	 & 	100	 & 	0.10	 & 	1	 & 	13	 & 	2.4	  	 & 	1.7	 	 & 	20.7	 & 	\textbf{5.3}	 & 	0.1\% 	 & 	1.4\% 	 & 	0.3\%\\\hline
w05\_100	 & 	10	 & 	100	 & 	0.50	 & 	9	 & 	199	 & 	78.0	 & 	37.4	 & 	371.7	 & 	\textbf{167.1}	 & 	0.9\% 	 & 	5.8\% 	 & 	3.2\%\\\hline
w09\_100	 & 	10	 & 	100	 & 	0.90	 & 	3	 & 	1243 & 	257.0	 & 	18.1	 & 	1114.1	 & 	\textbf{289.4}	 & 	0.1\% 	 & 	6.6\% 	 & 	3.7\%\\\hline
	\end{tabular}
	\caption{Numerical results obtained with sequential BiqBin on rudy instances for the Max-Cut problem.
		For an explanation of the columns see page~\pageref{tab1:main_text}.	}
	\label{tab:BiqBin_MC}
\end{table}

\begin{table}[h!]\scriptsize
\centering
	\begin{tabular}{|l||rrr|rrr|rrr|rrr|}\hline
\shortstack{instance \\ group} 	 & 	 \# inst. 	 & 	 $n$ 	 & 	 density 	 & 	 \shortstack{B\&B \\ (min)} 	 & 	 \shortstack{B\&B \\ (max)} 	 & 	 \shortstack{B\&B \\ (avg)} 	 & 	 \shortstack{time \\ (min)} 	 & 	 \shortstack{time \\ (max)} 	 & 	 \shortstack{time \\ (avg)} 	 & 	 \shortstack{init. gap \\ (min)} 	 & 	 \shortstack{init. gap \\ (max)} 	 & 	 \shortstack{init. gap \\ (avg)} \\\hline\hline
g05\_100	 & 	10	 & 	100	 & 	0.50	 & 	33	 & 	1751	 & 	367.0	 & 	44.3	 & 	1933.5	 & 	408.8	 & 	0.4\% 	 & 	1.2\% 	 & 	0.7\%\\\hline
pm1d\_100	 & 	10	 & 	100	 & 	0.99	 & 	39	 & 	1073	 & 	415.4	 & 	57.2	 & 	1302.3	 & 	519.6	 & 	3.2\% 	 & 	8.9\% 	 & 	5.7\%\\\hline
pm1s\_100	 & 	10	 & 	100	 & 	0.10	 & 	1	 & 	39	 & 	12.6	 & 	0.9	 & 	51.9	 & 	17.9	 & 	0.7\% 	 & 	4.0\% 	 & 	2.0\%\\\hline
pw01\_100	 & 	10	 & 	100	 & 	0.10	 & 	1	 & 	45	 & 	11.0	 & 	1.3	 & 	86.9	 & 	22.9	 & 	0.0\% 	 & 	0.9\% 	 & 	0.4\%\\\hline
pw05\_100	 & 	10	 & 	100	 & 	0.50	 & 	49	 & 	999	 & 	346.0	 & 	101.2	 & 	1155.5	 & 	453.0	 & 	0.4\% 	 & 	1.2\% 	 & 	0.8\%\\\hline
pw09\_100	 & 	10	 & 	100	 & 	0.90	 & 	113	 & 	589	 & 	280.8	 & 	170.4	 & 	824.2	 & 	399.2	 & 	0.4\% 	 & 	0.6\% 	 & 	0.5\%\\\hline
w01\_100	 & 	10	 & 	100	 & 	0.10	 & 	1	 & 	23	 & 	3.8	 & 	1.1	 & 	39.0	 & 	7.0	 & 	0.1\% 	 & 	2.5\% 	 & 	0.4\%\\\hline
w05\_100	 & 	10	 & 	100	 & 	0.50	 & 	31	 & 	401	 & 	210.4	 & 	58.9	 & 	564.0	 & 	304.2	 & 	2.3\% 	 & 	7.5\% 	 & 	4.7\%\\\hline
w09\_100	 & 	10	 & 	100	 & 	0.90	 & 	7	 & 	1667	 & 	423.2	 & 	17.3	 & 	2447.4	 & 	626.7	 & 	0.9\% 	 & 	8.1\% 	 & 	5.1\%\\\hline
	\end{tabular}
	\caption{Numerical results obtained with BiqCrunch on rudy instances for the Max-Cut problem.}
	\label{tab:BiqCrunch_MC}
\end{table}

\end{landscape}


\begin{landscape}

\begin{table}\scriptsize
\centering
	\begin{tabular}{|l||rrr|rrr|rrr|rrr|}\hline
\shortstack{instance \\ group} 	 & 	 \# inst. 	 & 	 $n$ 	 & 	 density 	 & 	 \shortstack{B\&B \\ (min)} 	 & 	 \shortstack{B\&B \\ (max)} 	 & 	 \shortstack{B\&B \\ (avg)} 	 & 	 \shortstack{time \\ (min)} 	 & 	 \shortstack{time \\ (max)} 	 & 	 \shortstack{time \\ (avg)} 	 & 	 \shortstack{init. gap \\ (min)} 	 & 	 \shortstack{init. gap \\ (max)} 	 & 	 \shortstack{init. gap \\ (avg)} \\\hline\hline
be100	 & 	10	 & 	100	 & 	1.00	 & 	1	 & 	9	 & 	2.2	 	 & 	7.4	 	 & 	81.7	 & 	26.2	 & 	0.0\% 	 & 	0.3\% 	 & 	0.1\%\\\hline
be120.3	 & 	10	 & 	120	 & 	0.30	 & 	1	 & 	1	 & 	1.0	 	 & 	5.6	 	 & 	31.5	 & 	11.7	 & 	0.0\% 	 & 	0.0\% 	 & 	0.0\%\\\hline
be120.8	 & 	10	 & 	120	 & 	0.80	 & 	1	 & 	35	 & 	10.4	 & 	17.9	 & 	267.4	 & 	\textbf{123.2}	 & 	0.0\% 	 & 	1.0\% 	 & 	0.3\%\\\hline
be150.3	 & 	10	 & 	150	 & 	0.30	 & 	1	 & 	191	 & 	30.8	 & 	16.3	 & 	832.9	 & 	\textbf{286.2}	 & 	0.0\% 	 & 	1.8\% 	 & 	0.4\%\\\hline
be150.8	 & 	10	 & 	150	 & 	0.80	 & 	5	 & 	223	 & 	68.8	 & 	221.8	 & 	1303.6	 & 	\textbf{656.6}	 & 	0.2\% 	 & 	1.2\% 	 & 	0.7\%\\\hline
be200.3	 & 	10	 & 	200	 & 	0.30	 & 	1	 & 	959	 & 	168.4	 & 	80.1	 & 	30241.9	 & 	\textbf{6085.5}	 & 	0.0\% 	 & 	2.3\% 	 & 	0.9\%\\\hline
be200.8	 & 	10	 & 	200	 & 	0.80	 & 	7	 & 	1095 & 	331.6	 & 	769.5	 & 	37449.9	 & 	\textbf{12280.0}	 & 	0.1\% 	 & 	2.2\% 	 & 	1.2\%\\\hline
be250	 & 	10	 & 	250	 & 	0.10	 & 	1	 & 	15	 & 	3.8	 	 & 	59.6	 & 	904.6	 & 	388.8	 & 	0.0\% 	 & 	0.3\% 	 & 	0.0\%\\\hline
bqp100	 & 	10	 & 	100	 & 	0.10	 & 	1	 & 	1	 & 	1.0	  	 & 	1.2	 	 & 	10.6	 & 	3.0	 & 	0.0\% 	 & 	0.0\% 	 & 	0.0\%\\ \hline
bqp250	 & 	10	 & 	250	 & 	0.10	 & 	1	 & 	81	 & 	9.8	 	 & 	81.7	 & 	9137.0	 & 	\textbf{1147.7}	 & 	0.0\% 	 & 	1.3\% 	 & 	0.1\%\\ \hline
	\end{tabular}
\caption{Numerical results obtained with sequential BiqBin on instances for unconstrained BQP.}\label{tab:BQP_biqBin}
\end{table}

\begin{table}\scriptsize
\centering
	\begin{tabular}{|l||rrr|rrr|rrr|rrr|}\hline
\shortstack{instance \\ group} 	 & 	 \# inst. 	 & 	 $n$ 	 & 	 density 	 & 	 \shortstack{B\&B \\ (min)} 	 & 	 \shortstack{B\&B \\ (max)} 	 & 	 \shortstack{B\&B \\ (avg)} 	 & 	 \shortstack{time \\ (min)} 	 & 	 \shortstack{time \\ (max)} 	 & 	 \shortstack{time \\ (avg)} 	 & 	 \shortstack{init. gap \\ (min)} 	 & 	 \shortstack{init. gap \\ (max)} 	 & 	 \shortstack{init. gap \\ (avg)} \\\hline\hline
be100	 & 	10	 & 	100	 & 	1.00	 & 	1	 & 	11	 & 	3.6	  	 & 	2.2	 	 & 	54.6	 & 	\textbf{17.2}	 & 	0.0\% 	 & 	0.7\% 	 & 	0.2\%\\\hline
be120.3	 & 	10	 & 	120	 & 	0.30	 & 	1	 & 	3	 & 	1.4	 	 &	1.8	 	 & 	26.3	 & 	\textbf{6.6}	 	 & 	0.0\% 	 & 	0.1\% 	 & 	0.0\%\\\hline
be120.8	 & 	10	 & 	120	 & 	0.80	 & 	3	 & 	57	 & 	17.2	 & 	15.2	 & 	422.1	 & 	126.0	 & 	0.0\% 	 & 	1.6\% 	 & 	0.7\%\\\hline
be150.3	 & 	10	 & 	150	 & 	0.30	 & 	1	 & 	79	 & 	25.6	 & 	4.5	 	 & 	841.0	 & 	295.4	 & 	0.0\% 	 & 	2.3\% 	 & 	0.7\%\\\hline
be150.8	 & 	10	 & 	150	 & 	0.80	 & 	19	 & 	141	 & 	62.8	 & 	226.2	 & 	1801.5	 & 	795.1	 & 	0.7\% 	 & 	1.7\% 	 & 	1.1\%\\\hline
be200.3	 & 	10	 & 	200	 & 	0.30	 & 	3	 & 	2253 & 	409.0	 & 	70.2	 & 	47534.1	 & 	9038.7	 & 	0.0\% 	 & 	2.8\% 	 & 	1.4\%\\\hline
be200.8	 & 	10	 & 	200	 & 	0.80	 & 	17	 & 	2803 & 	823.0	 & 	424.5	 & 	66334.4	 & 	20303.7	 & 	0.2\% 	 & 	2.7\% 	 & 	1.5\% \\\hline
be250	 & 	10	 & 	250	 & 	0.10	 & 	1	 & 	13	 & 	4.6	  	 & 	18.3	 & 	529.9	 & 	\textbf{185.7}	 & 	0.0\% 	 & 	0.4\% 	 & 	0.1\%\\\hline
bqp100	 & 	10	 & 	100	 & 	0.10	 & 	1	 & 	1	 & 	1.0	 	 & 	0.8	 	 & 	3.3	 	 & 	\textbf{1.3}	 	 & 	0.0\% 	 & 	0.0\% 	 & 	0.0\%\\\hline
bqp250	 & 	10	 & 	250	 & 	0.10	 & 	1	 & 	325	 & 	35.8	 & 	20.9	 & 	12812.7	 & 	1413.5	 & 	0.0\% 	 & 	1.9\% 	 & 	0.3\%\\\hline
	\end{tabular}
	\caption{Numerical results obtained with BiqCrunch on instances for unconstrained BQP.}
	\label{tab:BQP_BiqCrunch}
\end{table}

\end{landscape}


\begin{landscape}

\begin{table}\scriptsize
\centering
	\begin{tabular}{|l||rrr|rrr|rrr|rrr|}\hline
\shortstack{instance \\ group} 	 & 	 \# inst. 	 & 	 $n$ 	 & 	 density 	 & 	 \shortstack{B\&B \\ (min)} 	 & 	 \shortstack{B\&B \\ (max)} 	 & 	 \shortstack{B\&B \\ (avg)} 	 & 	 \shortstack{time \\ (min)} 	 & 	 \shortstack{time \\ (max)} 	 & 	 \shortstack{time \\ (avg)} 	 & 	 \shortstack{init. gap \\ (min)} 	 & 	 \shortstack{init. gap \\ (max)} 	 & 	 \shortstack{init. gap \\ (avg)} \\\hline\hline
120\_30\_0.25	 & 	5	 & 	120	 & 	0.25	 & 	19	 & 	57	 & 	30.6	 & 	105.9	 & 	196.6	 & 	143.7	 & 	1.9\% 	 & 	2.7\% 	 & 	2.2\%\\\hline
120\_30\_0.5	 & 	5	 & 	120	 & 	0.50	 & 	35	 & 	123	 & 	62.6	 & 	134.8	 & 	617.7	 & 	300.0	 & 	1.6\% 	 & 	2.1\% 	 & 	1.8\%\\\hline
120\_30\_0.75	 & 	5	 & 	120	 & 	0.75	 & 	13	 & 	289	 & 	111.4	 & 	47.6	 & 	1241.3	 & 	473.0	 & 	0.7\% 	 & 	1.9\% 	 & 	1.3\%\\\hline
120\_60\_0.25	 & 	5	 & 	120	 & 	0.25	 & 	1	 & 	25	 & 	11.8	 & 	13.7	 & 	113.2	 & 	62.6	 & 	0.2\% 	 & 	0.6\% 	 & 	0.4\%\\\hline
120\_60\_0.5	 & 	5	 & 	120	 & 	0.50	 & 	15	 & 	29	 & 	19.0	 & 	73.1	 & 	146.3	 & 	104.9	 & 	0.3\% 	 & 	0.4\% 	 & 	0.4\%\\\hline
120\_60\_0.75	 & 	5	 & 	120	 & 	0.75	 & 	1	 & 	13	 & 	7.8	 	 & 	8.5	 & 	77.8	 & 	38.0	 & 	0.1\% 	 & 	0.2\% 	 & 	0.2\%\\\hline
120\_90\_0.25	 & 	5	 & 	120	 & 	0.25	 & 	1	 & 	1	 & 	1.0	 & 	3.6	 & 	12.1	 & 	7.9	 & 	0.1\% 	 & 	0.1\% 	 & 	0.1\%\\\hline
120\_90\_0.5	 & 	5	 & 	120	 & 	0.50	 & 	3	 & 	19	 & 	9.8	 & 	38.5	 & 	184.5	 & 	103.7	 & 	0.0\% 	 & 	0.2\% 	 & 	0.1\%\\\hline
120\_90\_0.75	 & 	5	 & 	120	 & 	0.75	 & 	1	 & 	1	 & 	1.0	 & 	4.8	 & 	14.8	 & 	10.5	 & 	0.0\% 	 & 	0.0\% 	 & 	0.0\%\\\hline
140\_35\_0.25	 & 	5	 & 	140	 & 	0.25	 & 	25	 & 	261	 & 	125.4	 & 	134.1	 & 	1426.4	 & 	753.8	 & 	1.6\% 	 & 	2.9\% 	 & 	2.3\%\\\hline
140\_35\_0.5	 & 	5	 & 	140	 & 	0.50	 & 	83	 & 	1079	 & 	383.4	 & 	479.0	 & 	7274.1	 & 	2337.9	 & 	1.6\% 	 & 	2.7\% 	 & 	2.0\%\\\hline
140\_35\_0.75	 & 	5	 & 	140	 & 	0.75	 & 	159	 & 	1089	 & 	530.2	 & 	987.8	 & 	6486.0	 & 	3157.0	 & 	1.2\% 	 & 	1.7\% 	 & 	1.5\%\\\hline
140\_70\_0.25	 & 	5	 & 	140	 & 	0.25	 & 	3	 & 	131	 & 	45.4	 & 	38.0	 & 	634.3	 & 	245.1	 & 	0.3\% 	 & 	0.7\% 	 & 	0.5\%\\\hline
140\_70\_0.5	 & 	5	 & 	140	 & 	0.50	 & 	17	 & 	313	 & 	168.2	 & 	116.7	 & 	2098.5	 & 	1264.9	 & 	0.3\% 	 & 	0.6\% 	 & 	0.5\%\\\hline
140\_70\_0.75	 & 	5	 & 	140	 & 	0.75	 & 	3	 & 	49	 & 	19.4	 & 	35.5	 & 	153.4	 & 	87.3	 & 	0.1\% 	 & 	0.2\% 	 & 	0.2\%\\\hline
140\_105\_0.25	 & 	5	 & 	140	 & 	0.25	 & 	1	 & 	1	 & 	1.0	 & 	7.0	 & 	18.3	 & 	14.6	 & 	0.1\% 	 & 	0.1\% 	 & 	0.1\%\\\hline
140\_105\_0.5	 & 	5	 & 	140	 & 	0.50	 & 	1	 & 	23	 & 	7.0	 & 	16.5	 & 	274.7	 & 	95.0	 & 	0.0\% 	 & 	0.1\% 	 & 	0.1\%\\\hline
140\_105\_0.75	 & 	5	 & 	140	 & 	0.75	 & 	1	 & 	21	 & 	7.8	 & 	10.0	 & 	208.4	 & 	78.1	 & 	0.0\% 	 & 	0.1\% 	 & 	0.0\%\\\hline
160\_40\_0.25	 & 	5	 & 	160	 & 	0.25	 & 	49	 & 	773	 & 	298.6	 & 	421.4	 & 	8247.5	 & 	2885.2	 & 	1.6\% 	 & 	3.0\% 	 & 	2.3\%\\\hline
160\_40\_0.5	 & 	5	 & 	160	 & 	0.50	 & 	369	 & 	19307	 & 	4472.6	 & 	3510.2	 & 	170137.7	 & 	39019.6	 & 	1.6\% 	 & 	3.0\% 	 & 	2.1\%\\\hline
160\_40\_0.75	 & 	5	 & 	160	 & 	0.75	 & 	241	 & 	9353	 & 	3213.0	 & 	2600.9	 & 	71863.7	 & 	27203.1	 & 	1.0\% 	 & 	1.8\% 	 & 	1.5\%\\\hline
160\_80\_0.25	 & 	5	 & 	160	 & 	0.25	 & 	13	 & 	255	 & 	108.6	 & 	123.9	 & 	2091.8	 & 	885.8	 & 	0.4\% 	 & 	0.7\% 	 & 	0.5\%\\\hline
160\_80\_0.5	 & 	5	 & 	160	 & 	0.50	 & 	29	 & 	627	 & 	252.2	 & 	272.3	 & 	5303.6	 & 	\textbf{2080.1}	 & 	0.2\% 	 & 	0.5\% 	 & 	0.4\%\\\hline
160\_80\_0.75	 & 	5	 & 	160	 & 	0.75	 & 	3	 & 	3841	 & 	1045.0	 & 	44.9	 & 	27446.0	 & 	\textbf{7769.6}	 & 	0.1\% 	 & 	0.5\% 	 & 	0.3\%\\\hline

160\_120\_0.25	 & 	5	 & 	160	 & 	0.25	 & 	1	 & 	67	 & 	21.4	 & 	24.9	 & 	824.4	 & 	275.3	 & 	0.0\% 	 & 	0.2\% 	 & 	0.1\%\\\hline
160\_120\_0.5	 & 	5	 & 	160	 & 	0.50	 & 	1	 & 	51	 & 	18.6	 & 	38.7	 & 	611.2	 & 	248.8	 & 	0.0\% 	 & 	0.1\% 	 & 	0.1\%\\\hline
160\_120\_0.75	 & 	5	 & 	160	 & 	0.75	 & 	3	 & 	39	 & 	19.0	 & 	60.6	 & 	475.6	 & 	264.2	 & 	0.0\% 	 & 	0.1\% 	 & 	0.1\%\\\hline
	\end{tabular}
\caption{Numerical results obtained with sequential BiqBin on instances for the densest $k$-subgraph problem.}
\label{tab:densestBiqBin}
\end{table}

\begin{table}\scriptsize
\centering
	\begin{tabular}{|l||rrr|rrr|rrr|rrr|}\hline
\shortstack{instance \\ group} 	 & 	 \# inst. 	 & 	 $n$ 	 & 	 density 	 & 	 \shortstack{B\&B \\ (min)} 	 & 	 \shortstack{B\&B \\ (max)} 	 & 	 \shortstack{B\&B \\ (avg)} 	 & 	 \shortstack{time \\ (min)} 	 & 	 \shortstack{time \\ (max)} 	 & 	 \shortstack{time \\ (avg)} 	 & 	 \shortstack{init. gap \\ (min)} 	 & 	 \shortstack{init. gap \\ (max)} 	 & 	 \shortstack{init. gap \\ (avg)} \\\hline\hline
120\_30\_0.25	 & 	5	 & 	120	 & 	0.25	 & 	43	 & 	111	 & 	62.2	 & 	70.5	 & 	240.9	 & 	\textbf{132.0}	 & 	2.0\% 	 & 	2.9\% 	 & 	2.3\%\\\hline
120\_30\_0.5	 & 	5	 & 	120	 & 	0.50	 & 	55	 & 	185	 & 	108.6	 & 	82.5	 & 	332.3	 & 	\textbf{182.7}	 & 	1.7\% 	 & 	2.2\% 	 & 	2.0\%\\\hline
120\_30\_0.75	 & 	5	 & 	120	 & 	0.75	 & 	15	 & 	763	 & 	240.6	 & 	24.7	 & 	941.0	 & 	\textbf{312.3}	 & 	0.8\% 	 & 	2.0\% 	 & 	1.4\%\\\hline
120\_60\_0.25	 & 	5	 & 	120	 & 	0.25	 & 	3	 & 	61	 & 	27.0	 & 	7.2	 	 & 	128.2	 & 	\textbf{60.2}	 & 	0.3\% 	 & 	0.7\% 	 & 	0.5\%\\\hline
120\_60\_0.5	 & 	5	 & 	120	 & 	0.50	 & 	31	 & 	63	 & 	44.2	 & 	63.5	 & 	139.0	 & 	\textbf{99.5}	 & 	0.3\% 	 & 	0.4\% 	 & 	0.4\%\\\hline
120\_60\_0.75	 & 	5	 & 	120	 & 	0.75	 & 	7	 & 	31	 & 	19.4	 & 	12.3	 & 	81.3	 & 	\textbf{47.3}	 & 	0.1\% 	 & 	0.3\% 	 & 	0.2\%\\\hline
120\_90\_0.25	 & 	5	 & 	120	 & 	0.25	 & 	1	 & 	1	 & 	1.0	 & 	3.6	 		 & 	5.6	 	 & 	\textbf{4.4}	 & 	0.0\% 	 & 	0.1\% 	 & 	0.1\%\\\hline
120\_90\_0.5	 & 	5	 & 	120	 & 	0.50	 & 	3	 & 	15	 & 	7.8	 & 	11.6	 	 & 	69.9	 & 	\textbf{34.9}	 & 	0.0\% 	 & 	0.1\% 	 & 	0.1\%\\\hline
120\_90\_0.75	 & 	5	 & 	120	 & 	0.75	 & 	1	 & 	1	 & 	1.0	 	& 	2.0	 	& 	3.0		 & 	\textbf{2.5} 	 & 	0.0\% 	 & 	0.0\% 	 & 	0.0\%\\\hline
140\_35\_0.25	 & 	5	 & 	140	 & 	0.25	 & 	39	 & 	447	 & 	225.4	 & 	101.2	 & 	1137.6	 & 	\textbf{574.7}	 & 	1.8\% 	 & 	3.1\% 	 & 	2.5\%\\\hline
140\_35\_0.5	 & 	5	 & 	140	 & 	0.50	 & 	139	 & 	1813	 & 	613.8	& 308.1	 & 	3714.1	 & 	\textbf{1297.0}	 & 	1.8\% 	 & 	3.0\% 	 & 	2.2\%\\\hline
140\_35\_0.75	 & 	5	 & 	140	 & 	0.75	 & 	247	 & 	2063	 & 	985.0  & 610.8	 & 	3433.2	 & 	\textbf{1748.5}	 & 	1.3\% 	 & 	1.8\% 	 & 	1.6\%\\\hline
140\_70\_0.25	 & 	5	 & 	140	 & 	0.25	 & 	7	 & 	179	 & 	69.8	 & 	24.7	 & 	486.0	 & 	\textbf{200.6}	 & 	0.3\% 	 & 	0.8\% 	 & 	0.6\%\\\hline
140\_70\_0.5	 & 	5	 & 	140	 & 	0.50	 & 	37	 & 	877	 & 	415.0	 & 	108.7	 & 	2322.6	 & 	\textbf{1203.5}	 & 	0.3\% 	 & 	0.7\% 	 & 	0.5\%\\\hline
140\_70\_0.75	 & 	5	 & 	140	 & 	0.75	 & 	7	 & 	59	 & 	27.0	 & 	31.6	 & 	152.2	 & 	\textbf{76.3}	 & 	0.1\% 	 & 	0.2\% 	 & 	0.2\%\\\hline
140\_105\_0.25	 & 	5	 & 	140	 & 	0.25	 & 	1	 & 	1	 & 	1.0	 & 	4.1	 		& 	5.9		 & 	\textbf{5.0}	 & 	0.0\% 	 & 	0.1\% 	 & 	0.1\%\\\hline
140\_105\_0.5	 & 	5	 & 	140	 & 	0.50	 & 	1	 & 	11	 & 	4.2	 & 	4.4			 & 	67.5	 & 	\textbf{24.6}	 & 	0.0\% 	 & 	0.1\% 	 & 	0.0\%\\\hline
140\_105\_0.75	 & 	5	 & 	140	 & 	0.75	 & 	1	 & 	7	 & 	2.6	 & 	3.6			 & 	42.6	 & 	\textbf{15.1}	 & 	0.0\% 	 & 	0.1\% 	 & 	0.0\%\\\hline
160\_40\_0.25	 & 	5	 & 	160	 & 	0.25	 & 	73	 & 	1215	 & 	490.2	 & 256.3 & 	4317.5	 & 	\textbf{1697.0}	 & 	1.7\% 	 & 	3.1\% 	 & 	2.4\%\\\hline
160\_40\_0.5	 & 	5	 & 	160	 & 	0.50	 & 	509	 & 	25625	 & 	5825.8 & 1535.2	 & 	68554.5	 & 	\textbf{15918.8}	 & 	1.7\% 	 & 	3.0\% 	 & 	2.2\%\\\hline
160\_40\_0.75	 & 	5	 & 	160	 & 	0.75	 & 	327	 & 	12841	 & 	4636.6 & 1185.9	 & 	28700.9	 & 	\textbf{11046.9}	 & 	0.9\% 	 & 	1.9\% 	 & 	1.5\%\\\hline
160\_80\_0.25	 & 	5	 & 	160	 & 	0.25	 & 	31	 & 	485	 & 	213.4	 & 	126.1	 & 	1875.8	 & 	\textbf{817.2}	 & 	0.4\% 	 & 	0.8\% 	 & 	0.6\%\\\hline
160\_80\_0.5	 & 	5	 & 	160	 & 	0.50	 & 	33	 & 	1705	 & 	552.6	 & 	156.8	 & 	6592.5	 & 	2150.1	 & 	0.2\% 	 & 	0.6\% 	 & 	0.4\%\\\hline
160\_80\_0.75	 & 	5	 & 	160	 & 	0.75	 & 	3	 & 	8025	 & 	2191.8	 & 	17.6	 & 	29413.2	 & 	8229.8	 & 	0.0\% 	 & 	0.5\% 	 & 	0.3\%\\\hline
160\_120\_0.25	 & 	5	 & 	160	 & 	0.25	 & 	1	 & 	31	 & 	10.2	 & 	8.5		 & 	237.0	 & 	\textbf{81.3}	 & 	0.0\% 	 & 	0.2\% 	 & 	0.1\%\\\hline
160\_120\_0.5	 & 	5	 & 	160	 & 	0.50	 & 	1	 & 	21	 & 	7.0	 & 	7.1	 		 & 186.3	 & 	\textbf{62.4}	 & 	0.0\% 	 & 	0.1\% 	 & 	0.0\%\\\hline
160\_120\_0.75	 & 	5	 & 	160	 & 	0.75	 & 	1	 & 	9	 & 	4.6	 & 	7.3	  		 & 83.5		 & 	\textbf{40.2}	 & 	0.0\% 	 & 	0.0\% 	 & 	0.0\%\\\hline
	\end{tabular}
	\caption{Numerical results obtained with (tailored version of) BiqCrunch on instances for the densest $k$-subgraph problem.}
\label{tab:densestBiqCrunch}
\end{table}
\end{landscape}


\begin{landscape}

\begin{table}\scriptsize
\centering
	\begin{tabular}{|l||crr|c|rrr|c|rrr|c|rrr|}\hline
	&							&		&		&			\multicolumn{4}{c|}{BiqBin}	& \multicolumn{4}{c|}{Gurobi}	& \multicolumn{4}{c|}{SCIP}	\\
\shortstack{instance group} 	 & 	 \# inst. 	 & 	 $n$ 	 & 	 density 	 & 	\shortstack{solved \\ inst.} & \shortstack{time \\ (min)} 	 & 	 \shortstack{time \\ (max)} 	 & 	 \shortstack{time \\ (avg)} 	 & 	\shortstack{solved \\ inst.} & \shortstack{time \\ (min)} 	 & 	 \shortstack{time \\ (max)} 	 & 	 \shortstack{time \\ (avg)} 	 & 	\shortstack{solved \\ inst.} & \shortstack{time \\ (min)} 	 & 	 \shortstack{time \\ (max)} 	 & 	 \shortstack{time \\ (avg)} 	  \\\hline\hline
100\_A\_0\_1\_b\_10\_F\_-10\_10	 & 	15	 & 	100	 & 	0.95	 & 	\textbf{15}	 & 	26.3	 & 	5591.0	 & 	1085.8	 & 	2	 & 	1125.8	 & 	5715.4	 & 	3420.6	 & 	0	 & 	 -	 & 	 - 	 & 	 - 	\\\hline
100\_A\_0\_1\_b\_10\_F\_-5\_5	 & 	15	 & 	100	 & 	0.91	 & 	\textbf{14}	 & 	32.9	 & 	1376.4	 & 	406.5	 & 	1	 & 	5067.1	 & 	5067.1	 & 	5067.1	 & 	0	 & 	 -	 & 	 - 	 & 	 - 	\\\hline
100\_A\_0\_1\_b\_10\_F\_0\_10	 & 	15	 & 	100	 & 	0.91	 & 	\textbf{15}	 & 	1.9	 & 	3787.3	 & 	519.0	 	 & 	\textbf{15}	 & 	8.4	 & 	7279.1	 & 	951.7	 & 	2	 & 	1721.4	 & 	3805.2	 & 	2763.3	\\\hline
100\_A\_0\_1\_b\_10\_F\_0\_5	 & 	15	 & 	100	 & 	0.83	 & 	\textbf{15}	 & 	1.7	 & 	2090.3	 & 	321.9		 & 	14	 & 	11.4	 & 	2873.5	 & 	752.6	 & 	6	 & 	229.7	 & 	10363.3	 & 	2438.4	\\\hline
100\_A\_0\_1\_b\_15\_F\_-10\_10	 & 	15	 & 	100	 & 	0.95	 & 	\textbf{15}	 & 	19.5	 & 	5917.6	 & 	838.2	 & 	1	 & 	3187.6	 & 	3187.6	 & 	3187.6	 & 	0	 & 	 -	 & 	 - 	 & 	 - 	\\\hline
100\_A\_0\_1\_b\_15\_F\_-5\_5	 & 	15	 & 	100	 & 	0.91	 &  \textbf{15}	 & 	5.5	 & 	2016.9	 & 	630.6		 & 	0	 & 	 -	 & 	 - 	 & 	 - 	 & 	0	 & 	 -	 & 	 - 	 & 	 - 	\\\hline
100\_A\_0\_1\_b\_15\_F\_0\_10	 & 	15	 & 	100	 & 	0.91	 & 	\textbf{15}	 & 	20.2	 & 	2062.0	 & 	607.7	 & 	\textbf{15}	 & 	2.1	 & 	8325.4	 & 	738.4	 & 	4	 & 	3689.1	 & 	3789.7	 & 	3753.9	\\\hline
100\_A\_0\_1\_b\_15\_F\_0\_5	 & 	15	 & 	100	 & 	0.83	 & 	\textbf{15}	 & 	1.3	 & 	2360.2	 & 	331.9		 & 	\textbf{15}	 & 	2.5	 & 	9102.9	 & 	876.6	 & 	2	 & 	525.2	 & 	2048.2	 & 	1286.7	\\\hline
100\_A\_0\_1\_b\_20\_F\_-10\_10	 & 	15	 & 	100	 & 	0.95	 & 	\textbf{15}	 & 	8.1	 & 	1717.5	 & 	419.7		 & 	2	 & 	3520.0	 & 	8475.5	 & 	5997.8	 & 	0	 & 	 -	 & 	 - 	 & 	 - 	\\\hline
100\_A\_0\_1\_b\_20\_F\_-5\_5	 & 	15	 & 	100	 & 	0.91	 & 	\textbf{15}	 & 	13.7	 & 	8738.6	 & 	1078.1	 & 	0	 & 	 -	 & 	 - 	 & 	 - 	 & 	0	 & 	 -	 & 	 - 	 & 	 - 	\\\hline
100\_A\_0\_1\_b\_20\_F\_0\_10	 & 	15	 & 	100	 & 	0.91	 & 	\textbf{15}	 & 	3.3	 & 	2133.5	 & 	794.0		 & 	\textbf{15}	 & 	3.3	 & 	1761.8	 & 	155.4	 & 	1	 & 	7195.0	 & 	7195.0	 & 	7195.0	\\\hline
100\_A\_0\_1\_b\_20\_F\_0\_5	 & 	15	 & 	100	 & 	0.83	 & 	\textbf{15}	 & 	1.2	 & 	7835.7	 & 	830.8		 & 	\textbf{15}	 & 	4.6	 & 	946.3	 & 	125.3	 & 	3	 & 	3740.2	 & 	5295.7	 & 	4775.9	\\\hline
100\_A\_0\_3\_b\_10\_F\_-10\_10	 & 	15	 & 	100	 & 	0.95	 & 	7	 & 	37.1	 & 	9701.4	 & 	2479.7			 & 	\textbf{15}	 & 	19.0	 & 	2478.4	 & 	747.6	 & 	1	 & 	7814.7	 & 	7814.7	 & 	7814.7	\\\hline
100\_A\_0\_3\_b\_10\_F\_-5\_5	 & 	15	 & 	100	 & 	0.91	 & 	7	 & 	43.0	 & 	2449.8	 & 	1181.6			 & 	\textbf{15}	 & 	17.9	 & 	2677.0	 & 	555.3	 & 	6	 & 	6845.2	 & 	10224.6	 & 	8025.5	\\\hline
100\_A\_0\_3\_b\_10\_F\_0\_10	 & 	15	 & 	100	 & 	0.91	 & 	2	 & 	0.6	 & 	19.4	 & 	10.0				 & 	\textbf{15}	 & 	0.4	 & 	256.3	 & 	68.6	 & 	\textbf{15}	 & 	6.4	 & 	1478.6	 & 	619.7	\\\hline
100\_A\_0\_3\_b\_10\_F\_0\_5	 & 	15	 & 	100	 & 	0.83	 & 	4	 & 	1.3	 & 	132.0	 & 	34.7				 & 	\textbf{15}	 & 	0.3	 & 	271.8	 & 	62.5	 & 	\textbf{15}	 & 	9.9	 & 	2079.5	 & 	380.5	\\\hline
100\_A\_0\_3\_b\_15\_F\_-10\_10	 & 	15	 & 	100	 & 	0.95	 & 	6	 & 	89.5	 & 	8986.3	 & 	2388.7			 & 	\textbf{7}	 & 	2667.3	 & 	9917.8	 & 	7178.6	 & 	0	 & 	 -	 & 	 - 	 & 	 - 	\\\hline
100\_A\_0\_3\_b\_15\_F\_-5\_5	 & 	15	 & 	100	 & 	0.91	 & 	6	 & 	85.4	 & 	4443.2	 & 	1255.0			 & 	\textbf{8}	 & 	2491.6	 & 	9907.7	 & 	5893.9	 & 	0	 & 	 -	 & 	 - 	 & 	 - 	\\\hline
100\_A\_0\_3\_b\_15\_F\_0\_10	 & 	15	 & 	100	 & 	0.91	 & 	2	 & 	3.6	 & 	41.9	 & 	22.8				 & 	\textbf{15}	 & 	113.7	 & 	9264.4	 & 	3215.3	 & 	7	 & 	22.5	 & 	3710.7	 & 	983.6	\\\hline
100\_A\_0\_3\_b\_15\_F\_0\_5	 & 	15	 & 	100	 & 	0.83	 & 	4	 & 	10.9	 & 	69.6	 & 	34.5			 & 	\textbf{15}	 & 	1.8	 & 	10534.8	 & 	3634.6	 & 	6	 & 	12.7	 & 	10210.3	 & 	2831.7	\\\hline
100\_A\_0\_3\_b\_20\_F\_-10\_10	 & 	15	 & 	100	 & 	0.95	 & 	\textbf{5}	 & 	75.2	 & 	3770.8	 & 	1177.1	 & 	0	 & 	 -	 & 	 - 	 & 	 - 	 & 	0	 & 	 -	 & 	 - 	 & 	 - 	\\\hline
100\_A\_0\_3\_b\_20\_F\_-5\_5	 & 	15	 & 	100	 & 	0.91	 & 	\textbf{8}	 & 	85.6	 & 	10277.1	 & 	3027.1	 & 	1	 & 	5019.9	 & 	5019.9	 & 	5019.9	 & 	0	 & 	 -	 & 	 - 	 & 	 - 	\\\hline
100\_A\_0\_3\_b\_20\_F\_0\_10	 & 	15	 & 	100	 & 	0.91	 & 	3	 & 	10.5	 & 	26.9	 & 	16.8			 & 	\textbf{8}	 & 	201.2	 & 	5672.0	 & 	2115.6	 & 	5	 & 	108.5	 & 	3815.9	 & 	1668.9	\\\hline
100\_A\_0\_3\_b\_20\_F\_0\_5	 & 	15	 & 	100	 & 	0.83	 & 	5	 & 	0.9	 & 	422.2	 & 	129.4				 & 	\textbf{8}	 & 	161.0	 & 	8801.4	 & 	3224.8	 & 	6	 & 	198.2	 & 	2099.5	 & 	1370.9	\\\hline
100\_A\_-1\_1\_b\_0\_F\_-1\_1	 & 	15	 & 	100	 & 	 0.66 	 & 	\textbf{13}	 & 	7.4	 & 	6466.3	 & 	1643.1		 & 	3	 & 	2216.4	 & 	4496.0	 & 	3052.8	 & 	0	 & 	 -	 & 	 - 	 & 	 - 	\\\hline
100\_A\_-1\_1\_b\_0\_F\_-3\_3	 & 	15	 & 	100	 & 	 0.86 	 & 	\textbf{11}	 & 	15.3	 & 	7651.7	 & 	1909.6	 & 	2	 & 	1519.3	 & 	6271.7	 & 	3895.5	 & 	0	 & 	 -	 & 	 - 	 & 	 - 	\\\hline
100\_A\_-1\_1\_b\_0\_F\_-7\_7	 & 	15	 & 	100	 & 	 0.93 	 & 	\textbf{11}	 & 	21.1	 & 	9442.7	 & 	1652.8	 & 	6	 & 	5123.2	 & 	10320.5	 & 	7688.4	 & 	0	 & 	 -	 & 	 - 	 & 	 - 	\\\hline
100\_A\_-3\_3\_b\_0\_F\_-1\_1	 & 	15	 & 	100	 & 	 0.66 	 & 	\textbf{6}	 & 	6.6	 & 	9044.7	 & 	3187.7		 & 	2	 & 	940.7	 & 	5095.7	 & 	3018.2	 & 	0	 & 	 -	 & 	 - 	 & 	 - 	\\\hline
100\_A\_-3\_3\_b\_0\_F\_-3\_3	 & 	15	 & 	100	 & 	 0.86 	 & 	\textbf{6}	 & 	8.3	 & 	6538.8	 & 	2433.7		 & 	1	 & 	3958.7	 & 	3958.7	 & 	3958.7	 & 	0	 & 	 -	 & 	 - 	 & 	 - 	\\\hline
100\_A\_-3\_3\_b\_0\_F\_-7\_7	 & 	15	 & 	100	 & 	 0.93 	 & 	\textbf{4}	 & 	62.9	 & 	1753.4	 & 	674.8	 & 	2	 & 	9375.5	 & 	9945.0	 & 	9660.3	 & 	0	 & 	 -	 & 	 - 	 & 	 - 	\\\hline
100\_A\_-7\_7\_b\_0\_F\_-1\_1	 & 	15	 & 	100	 & 	 0.66 	 & 	\textbf{3}	 & 	220.3	 & 	3608.1	 & 	1710.5	 & 	2	 & 	2612.7	 & 	9769.2	 & 	6191.0	 & 	0	 & 	 -	 & 	 - 	 & 	 - 	\\\hline
100\_A\_-7\_7\_b\_0\_F\_-3\_3	 & 	15	 & 	100	 & 	 0.86 	 & 	\textbf{3}	 & 	35.1	 & 	1875.0	 & 	718.3	 & 	1	 & 	8122.2	 & 	8122.2	 & 	8122.2	 & 	0	 & 	 -	 & 	 - 	 & 	 - 	\\\hline
100\_A\_-7\_7\_b\_0\_F\_-7\_7	 & 	15	 & 	100	 & 	 0.93 	 & 	\textbf{3}	 & 	240.8	 & 	4196.8	 & 	1584.9	 & 	0	 & 	 -	 & 	 - 	 & 	 - 	 & 	0	 & 	 -	 & 	 - 	 & 	 - 	\\\hline
\multicolumn{1}{|r||}{total}	 & 	495	 & 		 & 		 & 	\textbf{298}	 & 		 & 		 & 		 & 	236	 & 		 & 		 & 		 & 	79	 & 		 & 		 & 		\\\hline
	\end{tabular}
	\caption{Comparison of sequential BiqBin, Gurobi and SCIP on randomly generated instances of \eqref{eqn:BQP}.
}
	\label{tab:RND_biqBin}
\end{table}

\end{landscape}

\subsection{Scaling properties of BiqBin}

In this section, we demonstrate how BiqBin is scaling across the high-performance computer that we used within the MPI framework.
We measured wall-times needed to solve each instance using the sequential solver and the parallel solver with 3, 6,\dots,48 CPU cores (i.e., MPI processes).
Note that one of the CPU cores was always reserved for the coordinator's tasks, while the others were used as workers.

In the ``sequential'' columns of Table~\ref{tab:scaling},
we report the computational times needed by the sequential algorithm and the times needed for the computations in the root node of the B\&B tree.
The former are used to compute the speed-up factor while the latter are considered as the times for the non-parallelizable parts of BiqBin
and are used to compute the upper bounds for the speed-up factors.

\begin{figure}[htp!]
	\centering
	\includegraphics[scale=0.6]{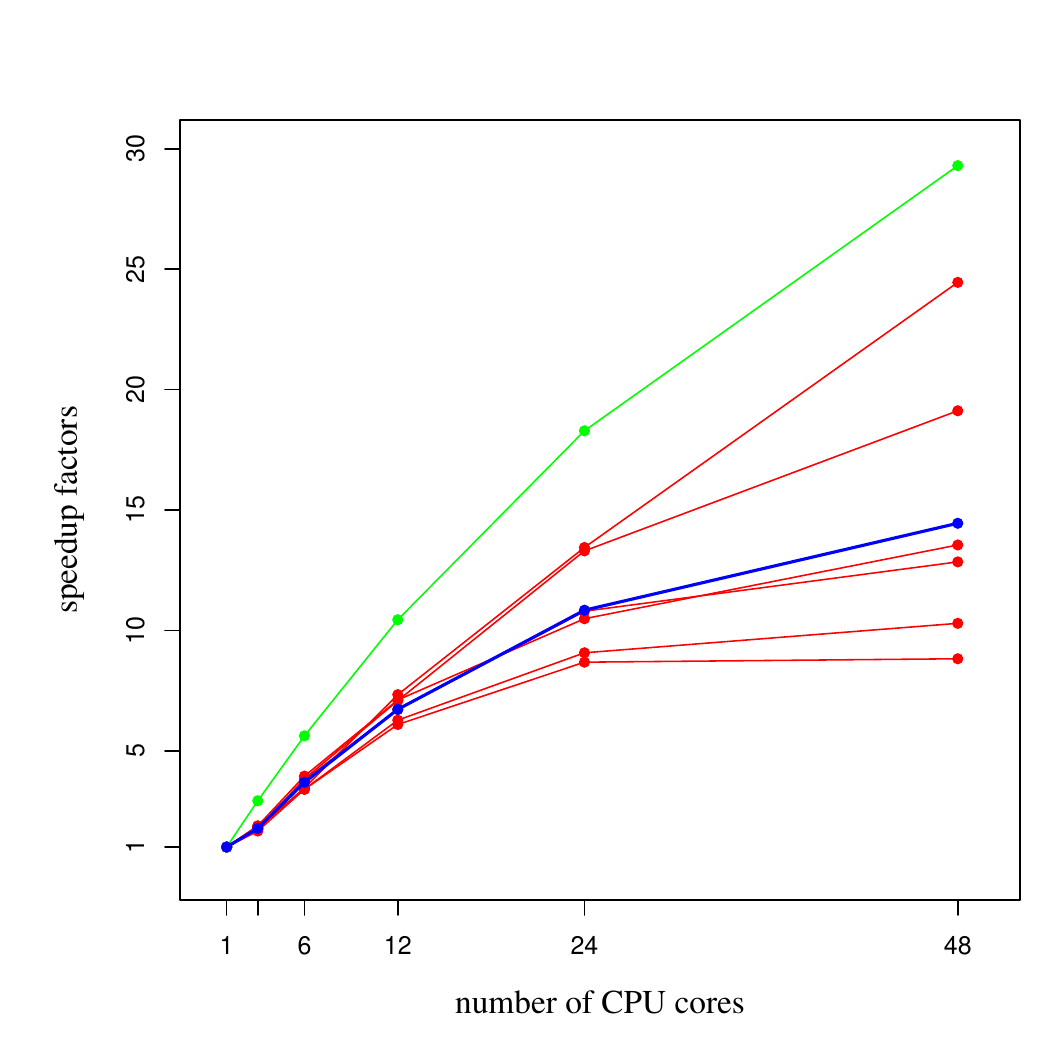}
	\caption{This plot shows how the code scales compared to upper-bounds $\mbox{UB}(n)$,
		represented by the green line.
		The red lines represent the scaling factors from Table~\ref{tab_agg_scaling_factors},
		while the blue line shows the average scaling factors (column-wise mean values of Table~\ref{tab_agg_scaling_factors}).
              }
	\label{fig:scaling_factors}	
\end{figure}

Pairs of columns denoted by ``3 cores'', ``6 cores'' etc. contain the sizes of the B\&B trees and the wall-times of the parallel BiqBin  using 3, 6,... CPU cores, respectively.
For each instance {\tt inst} from Table~\ref{tab:scaling},
we compute a vector of speed-up factors $\SU$, defined as (see, e.g.,~\cite[Eq. (3.4)]{navarro2014survey})
$$
\SU(i)=\frac{{\tt time}_1}{{\tt time}_i},
$$
where ${\tt time}_1$ is the time needed by the sequential solver while ${\tt time}_i$ denotes time needed by parallel solver using $i$ CPU cores.
These factors are aggregated in Table~\ref{tab_agg_scaling_factors} for each family of instances and are depicted in Figure~\ref{fig:scaling_factors}.

Times from the ``root time'' column  of Table~\ref{tab:scaling} are used to estimate the proportion of the BiqBin solver that is non-parallelizable.
We denote this estimate by $s$ and compute it as the sum of the 6th column  divided by the sum of the 5th column (result is $s=0.0136$).
The upper-bounds for the speed-up factors,  are computed according to Amdahl's law by formula (\cite[Eq. (3.6)]{navarro2014survey})
$$
	\mbox{UB}(n)=\frac{1}{s+(1-s)/n},
$$
where $n$ corresponds to the number of CPU cores  available.
These upper-bounds are depicted by the green curve in Figure~\ref{fig:scaling_factors}.
We can see that speed-up factors are close to the theoretical bound for the problems {g05} and  deviate a lot for the problems {pw09} and {w05}.

 Table~\ref{tab:scaling}  also contains the numbers of nodes in the B\&B trees generated by the parallel solver with different numbers of CPU cores.
 We can see that these numbers are different from the sequential solver and are varying with the number of CPU cores. The reason is mainly in the fact that parallel computations in different CPU cores are no more deterministic. Generating cutting planes
 in the process of solving \eqref{hyper_SDP}  includes random numbers. We fix the seed for the random number generator at the beginning of the computation in each CPU core. However, when we vary the number of CPU cores, the amount and the order of the computational work in the cores change. Consequently,
 different cutting planes might be applied to the same B\&B node and this results in slightly different bounds $\OPT_{\prob{HYP}}$ and finally in different sizes of B\&B trees.

\begin{landscape}	
		
		\begin{table}[htp!]\scriptsize
			\centering
			\begin{tabular}{|l||rr|rrr|rr|rr|rr|rr|rr|}\hline
				& 	  	 & 	 	 & 	  \multicolumn{3}{c|}{sequential} & \multicolumn{2}{c|}{3 cores} 	 & 	 \multicolumn{2}{c|}{6 cores} 	 & 	 \multicolumn{2}{c|}{12 cores} 	 & 	 \multicolumn{2}{c|}{24 cores} 	 & 	 \multicolumn{2}{c|}{48 cores} \\
				\shortstack{instance \\ name} 	 & 	 $n$ 	 & 	 density 	& \shortstack{B\&B}   &\shortstack{time}  & \shortstack{root~time} 	 &  \shortstack{B\&B} 	 & 	 \shortstack{time} 	 & 	 \shortstack{B\&B} 	 & 	 \shortstack{time} 	 & 	 \shortstack{B\&B} 	 & 	 \shortstack{time} 	 & 	 \shortstack{B\&B} 	 & 	 \shortstack{time} 	 & 	 \shortstack{B\&B} 	 & 	 \shortstack{time} 	 \\\hline\hline
				g05\_100.1 & 100 & 0.50 & 697 & 951.18 & 2.68 & 769 & 539.40 & 787 & 287.20 & 879 & 132.70 & 769 & 66.00 & 711 & 36.50 \\
				g05\_100.3 & 100 & 0.50 & 603 & 320.52 & 2.61 & 585 & 182.60 & 447 & 80.40 & 677 & 40.80 & 417 & 28.60 & 403 & 15.50  \\
				pm1d\_100.0 & 100 & 0.99 & 395 & 366.87 & 3.61 & 433 & 204.50 & 507 & 95.00 & 367 & 54.60 & 353 & 29.50 & 357 & 20.80  \\
				pm1d\_100.1 & 100 & 0.99 & 839 &  681.38 & 2.74 & 545 & 363.60 & 681 & 159.60 & 951 & 86.90 & 867 & 41.10 & 1045 & 31.70  \\
				pm1d\_100.2 & 100 & 0.99 & 551 & 372.83 & 3.70 & 493 & 220.80 & 621 & 105.70 & 443 & 56.40 & 525 & 32.10 & 369 & 16.70  \\
				pm1d\_100.4 & 100 & 0.99 & 557 & 357.24 & 3.72 & 439 & 228.60 & 445 & 104.20 & 571 & 53.20 & 435 & 31.00 & 371 & 23.80  \\
				pw05\_100.0 & 100 & 0.50 & 289 & 487.86 & 8.13 & 375 & 247.50 & 379 & 138.30 & 405 & 80.10 & 391 & 47.50 & 279 & 43.20  \\
				pw05\_100.6 & 100 & 0.50 & 263 & 503.85 & 6.71 & 183 & 279.80 & 183 & 113.00 & 215 & 59.40 & 313 & 47.10 & 279 & 30.00  \\
				pw09\_100.1 & 100 & 0.90 & 201 & 349.38 & 6.00 & 165 & 183.90 & 171 & 110.30 & 221 & 47.40 & 217 & 31.70 & 197 & 24.60  \\
				pw09\_100.5 & 100 & 0.90 & 137 &  260.52 & 6.65 & 107 & 181.20 & 133 & 66.40 & 133 & 51.70 & 121 & 38.10 &  99 & 43.60  \\
				pw09\_100.7 & 100 & 0.90 &  177 & 312.68 & 7.70 & 285 & 189.30 & 203 & 94.90 & 223 & 52.40 & 235 & 36.50 & 225 & 36.40 \\
				w05\_100.4 & 100 & 0.50 & 199 & 371.73 & 8.41 & 161 & 192.10 & 181 & 111.10 & 197 & 55.50 & 127 & 36.20 & 137 & 29.60  \\
				w05\_100.5 & 100 & 0.50 & 117 & 310.38 & 6.23 & 119 & 145.20 & 155 & 73.50 & 119 & 52.30 & 103 & 35.80 & 131 & 35.70 \\
				w05\_100.8 & 100 & 0.50 & 151 & 293.63 & 8.35 & 119 & 185.30 & 139 & 102.20 & 153 & 47.80 & 195 & 35.60 & 185 & 29.50  \\
				w09\_100.1 & 100 & 0.90 & 1243 & 1114.08 & 8.53 & 501 & 610.70 & 771 & 299.00 & 641 & 159.20 & 561 & 84.00 & 711 & 60.30 \\
				w09\_100.2 & 100 & 0.90 & 325 & 357.42 & 8.49 & 255 & 179.00 & 293 & 90.40 & 253 & 50.20 & 207 & 27.40 & 229 & 35.20  \\
				w09\_100.3 & 100 & 0.90 & 575 & 444.20 & 8.39 & 587 & 245.30 & 269 & 111.80 & 371 & 66.10 & 375 & 41.50 & 347 & 33.30  \\
				w09\_100.4 & 100 & 0.90 & 115 & 297.50 & 8.03 & 101 & 202.00 & 105 & 83.60 &  99 & 54.40 & 127 & 52.10 &  99 & 43.50 \\
				\hline
			\end{tabular}
			\caption{Numerical results obtained with parallel BiqBin for instances of the Max-Cut problem.
				}
			\label{tab:scaling}
		\end{table}

\begin{table}[htp!]\scriptsize
	\centering
	\begin{tabular}{|r||rrrrrrr|}
		\hline
		CPU cores & 1 & 3 & 6 & 12 & 24 & 48 & 96 \\ \hline\hline
g05 & 1.00 & 1.76 & 3.46 & 7.33 & 13.44 & 24.46 & 34.28 \\
pm1d & 1.00 & 1.75 & 3.83 & 7.08 & 13.30 & 19.12 & 24.70 \\
pw05 & 1.00 & 1.88 & 3.95 & 7.11 & 10.48 & 13.55 & 19.26 \\
pw09 & 1.00 & 1.66 & 3.40 & 6.09 & 8.68 & 8.82 & 11.87 \\
w05 & 1.00 & 1.87 & 3.40 & 6.27 & 9.07 & 10.29 & 12.48 \\
w09 & 1.00 & 1.79 & 3.78 & 6.71 & 10.80 & 12.85 & 17.59 \\
		\hline
	\end{tabular}
	\caption{This table contains aggregated scaling factors.
			The {g05} instances have on average best scaling factors, while instances {pw09}  and {w05} scale worst.}
		\label{tab_agg_scaling_factors}
\end{table}

\end{landscape}

\subsection{Solving large instances with parallel BiqBin}

In this section, we report numerical results obtained by parallel BiqBin on the Max-Cut, unconstrained BQP
and the densest $k$-subgraph instances that were to the best of our knowledge not solved so far.
We used 200 CPU cores and the results are collected in Tables~\ref{tab:largeMC}--\ref{tab:largeDkS}.
We can see that the hardest families of instances are again {g05} and {w09}.
To solve them, we need to compute more than one million B\&B nodes, which took at least 6~hours of wall-time.

\begin{table}[htp!]\scriptsize
\centering
	\begin{tabular}{|l||rrr|rr|}\hline
	 & 	 	 & 	 	 & 	 	 & 	 \multicolumn{2}{c|}{200 cores}\\
\shortstack{instance \\ group} 	 & 	 \# inst. 	 & 	 $n$ 	 & 	 density 	 & 	 \shortstack{B\&B \\ (avg)} 	 & 	 \shortstack{time \\ (avg)} \\\hline\hline
g05\_180	 & 	10	 & 	180	 & 	0.50	 & 	1712858.6	 & 	20700.5\\\hline
pm1d\_180	 & 	10	 & 	180	 & 	0.99	 & 	522354.4	 & 	10056.6\\\hline
pm1s\_180	 & 	10	 & 	180	 & 	0.10	 & 	60004.2	 & 	1163.1\\\hline
pw01\_180	 & 	10	 & 	180	 & 	0.10	 & 	4935.0	 & 	249.1\\\hline
pw05\_180	 & 	10	 & 	180	 & 	0.50	 & 	887735.0	 & 	20602.4\\\hline
pw09\_180	 & 	10	 & 	180	 & 	0.90	 & 	750231.8	 & 	18468.2\\\hline
w01\_180	 & 	10	 & 	180	 & 	0.10	 & 	36913.2	 & 	637.9\\\hline
w05\_180	 & 	10	 & 	180	 & 	0.50	 & 	674946.4	 & 	16216.1\\\hline
w09\_180	 & 	10	 & 	180	 & 	0.90	 & 	2061425.2	 & 	35320.7\\\hline
	\end{tabular}
	\caption{Numerical results obtained with parallel BiqBin for large instances of the Max-Cut problem.}\label{tab:largeMC}
\end{table}

\begin{table}[htp!]\scriptsize
\centering
	\begin{tabular}{|l||rrr|rr|}\hline
	 & 	 	 & 	 	 & 	 	 & 	 \multicolumn{2}{c|}{200 cores}\\
\shortstack{instance \\ group} 	 & 	 \# inst. 	 & 	 $n$ 	 & 	 density 	 & 	 \shortstack{B\&B \\ (avg)} 	 & 	 \shortstack{time \\ (avg)} \\\hline\hline
be250.3	 & 	10	 & 	250	 & 	0.30	 & 	1367.4	 & 	1133.0\\\hline
be250.8	 & 	10	 & 	250	 & 	0.80	 & 	7838.8	 & 	3583.2\\\hline
be300.3	 & 	10	 & 	300	 & 	0.30	 & 	5462.8	 & 	4207.5\\\hline
be300.8	 & 	10	 & 	300	 & 	0.80	 & 	110806.8	 & 	74765.2\\\hline
	\end{tabular}
	\caption{Numerical results obtained with parallel BiqBin for large instances of the unconstrained BQP.}\label{tab:BiqBQP}
\end{table}

\begin{table}[htp!]\scriptsize
	\centering
	\begin{tabular}{|l||rrr|rr|}\hline
	 & 	 	 & 	 	 & 	 	 & 	 \multicolumn{2}{c|}{200 cores}\\
\shortstack{instance \\ group} 	 & 	 \# inst. 	 & 	 $n$ 	 & 	 density 	 & 	 \shortstack{B\&B \\ (avg)} 	 & 	 \shortstack{time \\ (avg)} \\\hline\hline
120\_30\_0.25	 & 	5	 & 	120	 & 	0.25	 & 	33.8	 & 	32.1\\\hline
120\_30\_0.5	 & 	5	 & 	120	 & 	0.50	 & 	70.2	 & 	33.9\\\hline
120\_30\_0.75	 & 	5	 & 	120	 & 	0.75	 & 	117.0	 & 	29.0\\\hline
120\_60\_0.25	 & 	5	 & 	120	 & 	0.25	 & 	17.4	 & 	27.8\\\hline
120\_60\_0.5	 & 	5	 & 	120	 & 	0.50	 & 	21.4	 & 	34.7\\\hline
120\_60\_0.75	 & 	5	 & 	120	 & 	0.75	 & 	10.2	 & 	20.2\\\hline
120\_90\_0.25	 & 	5	 & 	120	 & 	0.25	 & 	1.0	 & 	9.8\\\hline
120\_90\_0.5	 & 	5	 & 	120	 & 	0.50	 & 	12.6	 & 	45.0\\\hline
120\_90\_0.75	 & 	5	 & 	120	 & 	0.75	 & 	1.0	 & 	12.5\\\hline
140\_35\_0.25	 & 	5	 & 	140	 & 	0.25	 & 	146.2	 & 	51.5\\\hline
140\_35\_0.5	 & 	5	 & 	140	 & 	0.50	 & 	395.8	 & 	59.7\\\hline
140\_35\_0.75	 & 	5	 & 	140	 & 	0.75	 & 	479.4	 & 	59.4\\\hline
140\_70\_0.25	 & 	5	 & 	140	 & 	0.25	 & 	58.2	 & 	37.9\\\hline
140\_70\_0.5	 & 	5	 & 	140	 & 	0.50	 & 	173.8	 & 	56.0\\\hline
140\_70\_0.75	 & 	5	 & 	140	 & 	0.75	 & 	12.6	 & 	29.3\\\hline
140\_105\_0.25	 & 	5	 & 	140	 & 	0.25	 & 	1.0	 & 	17.9\\\hline
140\_105\_0.5	 & 	5	 & 	140	 & 	0.50	 & 	9.8	 & 	50.3\\\hline
140\_105\_0.75	 & 	5	 & 	140	 & 	0.75	 & 	7.4	 & 	37.6\\\hline
160\_40\_0.25	 & 	5	 & 	160	 & 	0.25	 & 	395.0	 & 	84.5\\\hline
160\_40\_0.5	 & 	5	 & 	160	 & 	0.50	 & 	4817.8	 & 	329.9\\\hline
160\_40\_0.75	 & 	5	 & 	160	 & 	0.75	 & 	3398.6	 & 	234.5\\\hline
160\_80\_0.25	 & 	5	 & 	160	 & 	0.25	 & 	113.4	 & 	62.1\\\hline
160\_80\_0.5	 & 	5	 & 	160	 & 	0.50	 & 	253.8	 & 	76.8\\\hline
160\_80\_0.75	 & 	5	 & 	160	 & 	0.75	 & 	1071.4	 & 	92.1\\\hline
160\_120\_0.25	 & 	5	 & 	160	 & 	0.25	 & 	25.4	 & 	84.6\\\hline
160\_120\_0.5	 & 	5	 & 	160	 & 	0.50	 & 	29.4	 & 	72.7\\\hline
160\_120\_0.75	 & 	5	 & 	160	 & 	0.75	 & 	17.8	 & 	76.5\\\hline
180\_45\_0.25	 & 	5	 & 	180	 & 	0.25	 & 	1548.2	 & 	188.6\\\hline
180\_45\_0.5	 & 	5	 & 	180	 & 	0.50	 & 	2109.4	 & 	265.0\\\hline
180\_45\_0.75	 & 	5	 & 	180	 & 	0.75	 & 	29243.4	 & 	2089.7\\\hline
180\_90\_0.25	 & 	5	 & 	180	 & 	0.25	 & 	923.0	 & 	144.9\\\hline
180\_90\_0.5	 & 	5	 & 	180	 & 	0.50	 & 	1167.4	 & 	165.8\\\hline
180\_90\_0.75	 & 	5	 & 	180	 & 	0.75	 & 	2623.8	 & 	246.9\\\hline
180\_135\_0.25	 & 	5	 & 	180	 & 	0.25	 & 	11.8	 & 	79.5\\\hline
180\_135\_0.5	 & 	5	 & 	180	 & 	0.50	 & 	29.4	 & 	133.2\\\hline
180\_135\_0.75	 & 	5	 & 	180	 & 	0.75	 & 	29.4	 & 	113.2\\\hline
200\_50\_0.25	 & 	5	 & 	200	 & 	0.25	 & 	7400.6	 & 	849.4\\\hline
200\_50\_0.5	 & 	5	 & 	200	 & 	0.50	 & 	24145.4	 & 	2578.5\\\hline
200\_50\_0.75	 & 	5	 & 	200	 & 	0.75	 & 	53016.2	 & 	4980.4\\\hline
200\_100\_0.25	 & 	5	 & 	200	 & 	0.25	 & 	1392.2	 & 	254.9\\\hline
200\_100\_0.5	 & 	5	 & 	200	 & 	0.50	 & 	3801.8	 & 	501.3\\\hline
200\_100\_0.75	 & 	5	 & 	200	 & 	0.75	 & 	3067.4	 & 	358.6\\\hline
200\_150\_0.25	 & 	5	 & 	200	 & 	0.25	 & 	43.4	 & 	145.1\\\hline
200\_150\_0.5	 & 	5	 & 	200	 & 	0.50	 & 	129.0	 & 	177.1\\\hline
200\_150\_0.75	 & 	5	 & 	200	 & 	0.75	 & 	621.0	 & 	260.3\\\hline
	\end{tabular}
	\caption{Numerical results obtained with parallel BiqBin for large instances of the densest $k$-subgraph problem.}\label{tab:largeDkS}
\end{table}

\section{Conclusions}

In this paper we describe BiqBin, a solver for linearly constrained quadratic problems,
which is capable to solve to optimality instances that are due to their size unsolvable by other existing methods and tools.
The main idea underlying this solver is the exact penalty reformulation of a \eqref{eqn:BQP} instance to an instance of \eqref{eqn:MC},
introduced by Lasserre and enhanced by two co-authors of this paper in~\cite{GuWi:19}.

We provide necessary theoretical results needed to explain the work-flow of the problem reformulations,
relaxations, and finally the details related to the C implementation of BiqBin as an efficient parallel solver.
The solver is also available as a web service, which is connected to the \rew{high-performance} computer at the University of Ljubljana, Faculty of mechanical engineering.

We present extensive numerical results, where BiqBin is benchmarked against BiqCrunch, GUROBI and SCIP.
It can be concluded that BiqBin is outperforming other solvers on the Max-Cut instances,
that it is competitive with BiqCrunch on the instances of unconstrained binary quadratic problems,
and that it is slightly worse than BiqCrunch on the instances of the densest $k$-subgraph problem. The latter is
expected since BiqCrunch is specially adapted to solve problems of this type.
On these three families of instances GUROBI and SCIP are non-competitive.

However, when the number of linear constraints slightly increases,
like it happens on the fourth family of benchmark instances (randomly generated instances of \eqref{eqn:BQP}),
GUROBI and SCIP become solvers of interest.
They solve the problems in the original formulation and the initial linear constraints are important in generating new cutting planes. Therefore, a larger number of linear constraints usually results in a better performance of these solvers, while BiqBin reformulates the problem into an instance of Max-Cut, hence the feasible set always consists of all possible binary vectors.
However, on the benchmark instances, BiqBin is still best-performing, while BiqCrunch demonstrates very weak performance and was eliminated from the reported numerical results.

We showed that the BiqBin solver scales very well,
hence \ifAcm \else \linebreak \fi high-performance computers are the infrastructure to be used to solve instances,
that are out of reach for other (sequential) state-of-the-art solvers.

As part of our future work,
we plan to merge BiqBin, which uses the bundle method as the computational core, with MADAM introduced in~\cite{HrgaPovh:20},
where an alternating direction method of multipliers is developed for solving hard semidefinite relaxations.
Additionally, the MPI communication among the processes can be further simplified to enable efficient scaling over larger HPC systems.
This can be achieved by employing a one-sided MPI communication. We will also further improve the performance of BiqBin by enhancing  an early  stopping condition for the case of  infeasibility of \eqref{eqn:BQP}.

However, to go much further and solve much larger instances of \eqref{eqn:BQP} or even more general classes of discrete optimization problems,
we need to enhance the existing approach with new advances from polynomial optimization, artificial intelligence and problem specific theoretical findings.

	\section*{Acknowledgements}
        This project was supported by the Austrian Science Fund (FWF):
        I\,3199-N31,
		and by the Slovenian Research Agency: Program P2-0162, and the projects N1-0057, N1-0071, J1-2453, and J1-1691.
		The second author also acknowledges the financial support from the Slovenian Research Agency under the program for young researchers (MR+).
		Furthermore, this project has received funding from
		the European Union's Horizon~2020 research and
		innovation programme under the
		Marie Sk\l{}odowska-Curie grant agreement No~764759.
		Fruitful discussions with Leon Kos, Franz Rendl, and Alen Vegi Kalamar helped a lot to resolve several challenges during the project.
                \rew{Finally, we thank two anonymous referees for improving an earlier version of this work.}
\ifAcm
\else
	\bibliography{BiqBin_biblio}
	\bibliographystyle{plain}
\fi
\end{document}